\documentclass[default]{svmult}

\usepackage{mathptmx}       
\usepackage{helvet}         
\usepackage{courier}        
\usepackage{type1cm}        

\usepackage{makeidx}         
\usepackage{graphicx}        
\usepackage{multicol}        

\usepackage{amsmath,amssymb,amsfonts}
\usepackage{wasysym} 
\usepackage{cancel}

\newcommand{\dive}{\textrm{div}}
\newcommand{\grad}{\textrm{grad}}
\newcommand{\Dive}{\textrm{Div}}
\newcommand{\Grad}{\textrm{Grad}}

\newtheorem{thm}{Theorem}[section]
\newtheorem{defn}{Definition}[section]
\newtheorem{prop}{Proposition}[section]

\makeindex     

\begin{document}
\baselineskip=13pt

\title*{Analysis and numerics of the propagation speed\\ for hyperbolic reaction-diffusion models}
\titlerunning{Analysis and numerics for hyperbolic reaction-diffusion models}
\author{Corrado Lattanzio, Corrado Mascia, Ramon G. Plaza, Chiara Simeoni}
\authorrunning{C. Lattanzio, C. Mascia, R.G. Plaza, C. Simeoni}
\institute{Corrado Lattanzio \at Dipartimento di Ingegneria e Scienze dell'Informazione e Matematica, Universit\`a degli Studi dell'Aquila (Italy), \email{corrado@univaq.i}
\and Corrado Mascia \at Dipartimento di Matematica ``Guido Castelnuovo", Sapienza Universit\`a di Roma (Italy), \email{corrado.mascia@uniroma1.it}
\and Ramon G. Plaza \at Instituto de Investigaciones en Matem\'aticas Aplicadas y en Sistemas, Universidad Nacional Aut\'onoma de M\'exico (M\'exico), \email{plaza@mym.iimas.unam.mx}
\and Chiara Simeoni \at Laboratoire J.A. Dieudonn\'e, Universit\'e C\^ote d'Azur, Nice (France), \email{chiara.simeoni@univ-cotedazur.fr}}

\maketitle

\abstract{
In this paper, we analzye propagating fronts in the context of hyperbolic theories of dissipative processes.
These can be considered as a natural alternative to the more classical parabolic models.
Emphasis is given toward the numerical computation of the invasion velocity.
The first Section is devoted to the presentation of different models for reaction-diffusion phenomena,
supporting the idea of the advantages of a description based on hyperbolic equations.
Among other advantages, such modeling could provide a detailed description of the transient
dynamics of the phenomenon under observation.
Three basic numerical schemes are also presented; two of them can, in principle, be applied
to general hyperbolic systems, at the price of reduced performances when dealing
with discontinuous initial data.
In the second Section, we focus on a specific class of $2\times 2$ system corresponding to second
order partial differential equations in one space dimension, adapted for simplified modeling
of reaction-diffusion equations.
Specifically, we focus on notable traveling wave solutions, called {\it propagation fronts}.
Particular cases where the speed of propagation can be explicitly computed are also provided.
The third (and final) Section starts with the presentation of the {\it phase-plane algorithm}
which bears a reliable approximation of the propagation speed, assessing its validity in the case 
with damping where an explicit formula is available.
Then, we propose two PDE-based algorithms to approximate such velocity, named,
respectively, {\it scout\&spot algorithm} (based on tracking the level curve of some intermediate value of the profile)
and {\it LeVeque--Yee formula} (given by the average value of the discrete transport velocity).
Finally, we attest the well-foundedness of both the approaches and conclude by suggesting the second one as more
efficient tool in the determination of the speed.
}



\section{Models for reaction-diffusion phenomena}
\label{sect:models}

\setcounter{equation}{0}
\renewcommand{\theequation}{1.\arabic{equation}}

In this Section, we present different type of models useful for describing reaction-diffusion phenomena.
The standard approach gives raise to a parabolic equation which is very well suited to explain simple
events such as heat transmission in close-to-equilibrium regime.
In the standard linear case, such modeling has been criticised for three main reasons:
\begin{itemize}
\item[1.] infinite speed of propagation;
\item[2.] lack of time-delay and related inertial effects;
\item[3.] excepionality of well-posed boundary value problems.
\end{itemize}
In addition to the discussion relative to inertia (started by Eckart
in the 40s \cite{Ecka40} and continued in \cite{GeroLind90,LehnReulRubi18} in the context of relativity), 
other fields where the hyperbolic terms are relevant for applications are, among others,
in biological tissues \cite{DunbOthm86, OthmDunbAlt88, XuLu11, XuSeffLu08}, population growth \cite{MendCama97}, 
forest fire models \cite{MendLleb97}...

Here, starting from Subsection \ref{subsec:scalar} (dealing with scalar equations) and proceeding with Subsection
\ref{subsec:systems} (focusing on systems), we follow the point of view that a description making use of hyperbolic equations
--starting from the basic example of the {\it telegraph equation}-- is viable and more appropriate when the relaxation time required 
to sense the change of the overall phenomenon is sufficiently large as compared to the diffusivity coefficient.
Indeed, differences may emerge in the transient time, whose cumulation may influence significantly the final outcome.

Section ends with a presentation of three different numerical schemes which can be easily implemented
in order to obtain reliable approximation of a reaction-diffusion model of hyperbolic type.
We stress that we do not regard hyperbolic numerical schemes as a tool for approximating parabolic equations;
rather, we focus on hyperbolic models considered as a different language useful for describing dissipative mechanisms in a
modified manner which could be interesting in the modelling of distinct phenomena in far-from-equilibrium regimes.

\subsection{Diffusion is not always a parabolic mechanism}
\label{subsec:scalar}

The standard approach to heat conduction in a homogeneous medium is based on the continuity relation linking
the scalar unknown variable $u$ with the vector-valued flux function $\mathbf{v}$, by means of the balance identity
\begin{equation*}\label{integralidentity}
	\frac{d}{dt}\int_{\Omega} u(\mathbf{x},t)\,d\mathbf{x}+\int_{\partial\Omega} \mathbf{v}\cdot\mathbf{n}\,d\sigma
		=\int_{\Omega} f \,d\mathbf{x},
\end{equation*}
where $\Omega$ is an arbitrarily chosen control region with $d\mathbf{x}$ corresponding volume element,
$\mathbf{n}$ is the outward normal to the smooth boundary $\partial\Omega$ with $d\sigma$ boundary element,
and $f$ is a volume contribution, to be considered, at first, as a given external constraint.

Applying Divergence Theorem, we can consider the localised version
\begin{equation}\label{continuity}
	\partial_t u + \dive_{\mathbf{x}}\, \mathbf{v} = f,
\end{equation}
where $u$ and $\mathbf{v}$ describe respectively (heat) density and (heat) flux.
The former is a scalar quantity; the latter is a vector with same dimension of the space variable $\mathbf{x}$.

To provide a closed system, equation \eqref{continuity} has to be coupled with some relation
between $u$ and $\mathbf{v}$.
A frequent choice is the {\it Fourier's law}
\begin{equation}\label{fourier}
	\mathbf{v}=-a\,\grad_{\mathbf{x}} u
\end{equation}
for some non-negative proportionality parameter $a$,
which may explicitly depend on space $\mathbf{x}$ and time $t$ --as in the case of heterogeneous media--
and also on the density variable itself $u$ and its derivatives.
Here, we focus mainly on the case where $a$ is a given positive constant, i.e. $a>0$.
Linear relation \eqref{fourier} is also called {\it Fick's law} when considered in bio-mathematical settings,
{\it Ohm's law} in electromagnetism, and  {\it Darcy's law} in porous media.

Coupling identity \eqref{continuity} with relation \eqref{fourier} gives raise to the balance law
\begin{equation}\label{scalar_par_rd}
	\partial_t u =\dive_{\mathbf{x}}\left(a\,\grad_{\mathbf{x}} u\right)+f.
\end{equation}
While the continuity equation \eqref{continuity} can be considered reliable in general contexts,
equation \eqref{fourier} should be regarded as a single possible choice among many others.
In fact, quoting Lars Onsager (see \cite{Onsa31}), 
{\it Fourier's law is an approximate description of the process
of conduction, which neglects the (short) relaxation time $\tau$ needed for acceleration}.
For practical purposes (as in heat conduction) the time-lag can be neglected in all cases that are likely to be studied.
Nevertheless, in many applications --among others, for far-from-equilbrium regimes, such as the study of living tissues
and thermal resonance-- extensions of the Fourier's law are required, with the specific aim of providing a more robust model.

A first significant alternative to \eqref{fourier} is supported by the intuition that a delayed version should hold in place
of the instantaneous response.
The fact that the system requires a strictly positive amount of time $\tau$ to sense the gradient change
translates into an identity of the {\it phase-lag relationship}
\begin{equation}\label{delay}
	\mathbf{v} (\mathbf{x},t+\tau) = - a\,\grad_{\mathbf{x}} u(\mathbf{x},t).
\end{equation}
Unfortunately, as proved in \cite{JordDaiMick08}, the phase-lag model is ill-posed in the sense of Hadamard
since it lacks of continuous dependence with respect to the initial data (see also \cite{DrehQuinRack09}).

Surprisingly enough, well-posedness can be restored by truncating the Taylor's expansion for the unknown $\mathbf{v}$.
Assuming $\tau$ to be small, we can consider the approximation
\begin{equation*}\label{taylorv_scalar}
	\begin{aligned}
	\mathbf{v} (\mathbf{x},t+\tau)
		&= \mathbf{v} (\mathbf{x},t)+\tau \partial_t \mathbf{v} (\mathbf{x},t)+o(\tau)\\
		& \approx  \mathbf{v} (\mathbf{x},t)+\tau\,\partial_t \mathbf{v} (\mathbf{x},t),
	\end{aligned}
\end{equation*}
giving raise to the {\it Maxwell--Cattaneo's law}.
Putting together with the balance law \eqref{continuity}, we obtain the {\it (hyperbolic) reaction-diffusion system with relaxation}
\begin{equation}\label{mc_scalar}
	\left\{\begin{aligned}
	\partial_t u + \dive_{\mathbf{x}}\mathbf{v} & = f,\\
	\tau\partial_t \mathbf{v} + a\,\grad_{\mathbf{x}} u & = - \mathbf{v}.
	\end{aligned}\right.
\end{equation}
The Maxwell--Cattaneo's law can be considered as a way for incorporating into the diffusion modelling some additional
physical terms arising in the framework of Extended Irreversible Thermodynamics, \cite{CimmJouRuggVan14,JouCasaVazq10}.
Such law, to be considered as a constitutive identity, has been originally proposed by Cattaneo \cite{Catt48,Catt58},
following some pioneering intuition of James Clerk Maxwell (among others, let us quote \cite{JosePrez89,MorsFesh53}).
Sometimes, equation \eqref{mc_scalar} is attributed to Vernotte \cite{Vern58}, and --more rarely-- to Chester \cite{Ches63}.
Extensions has been also proposed in \cite{Chri09}.

Eliminating the unknown $\mathbf{v}$ in the coupled system  \eqref{continuity} and \eqref{mc_scalar}, 
we obtain the {\it one-field equation}, namely
\begin{equation}\label{scalar_hhd}
	\tau \partial_{tt} u + \partial_t\left(u-\tau f\right) = \dive_{\mathbf{x}}\left(a\,\grad_{\mathbf{x}} u\right)+f.
\end{equation}
The focal idea is that the balance between the flux $\mathbf{v}$ and the gradient $\grad_{\mathbf{x}} u$ of the density $u$
is achieved only asymptotically in time, with decay described by the {\it relaxation time} $\tau>0$.
Such quantity can be regarded as the characteristic time for the crossover between ballistic motion and the onset of diffusion. 

The Maxwell--Cattaneo's law furnishes the differential version of the delayed response to a change in the gradient $\grad_{\mathbf{x}} u$
as described by a memory kernel given by the exponential-rate law
\begin{equation*}\label{memory0}
	\mathbf{v}(\mathbf{x},t)=\mathbf{v}_0(\mathbf{x})e^{-t/\tau}-\frac{1}{\tau}\int_0^t e^{-(t-s)/\tau} a\,\grad_{\mathbf{x}} u(\mathbf{x},s)\,ds
\end{equation*}
which corresponds to the analogous formula in the context of viscoelasticity.
Incidentally, let us observe that the nonlocality of the time-integral --to be compared with the instantaneous
relationship \eqref{delay}-- can be regarded as a partial justification of the fact that the reaction-diffusion system \eqref{mc_scalar}
is proved to be time-locally well-posed.

The main flaw is that equation \eqref{scalar_hhd} can violate the second law of thermodynamics, admitting scenarios where heat appear
to be moving from cold to hot (see \cite{KornBerg98}).
In this respect, correction to the notion of entropy have been proposed in order to partially solve the problem
(for the case with no source term, see \cite{CriaLleb93}). 

An alternative approach is based on the postulation that the usual continuity equation \eqref{continuity} should be replaced by a delayed identity
\begin{equation*}
	\partial_t u (\mathbf{x},t+\tau) + \dive_{\mathbf{x}}\, \mathbf{v} (\mathbf{x},t)= f(\mathbf{x},t).
\end{equation*}
Truncating again the Taylor's expansion for $u$ with respect to the second argument, we end up with
\begin{equation}\label{delayed_continuity}
	\tau\,\partial_{tt} u +\partial_t u + \dive_{\mathbf{x}}\, \mathbf{v} = f.
\end{equation}
Then, coupling with the standard Fourier's law \eqref{fourier}, equation \eqref{delayed_continuity}
gives the so-called {\it (hyperbolic) reaction-diffusion equation with damping}
\begin{equation}\label{damp_rd}
	\tau\,\partial_{tt} u+\partial_t u
		= \dive_{\mathbf{x}}\left(a\,\grad_{\mathbf{x}} u\right)+f.
\end{equation}
An alternative approach leading to a variation of \eqref{damp_rd} is proposed in \cite{AliZhan05}, 
where the hyperbolic equation \eqref{damp_rd} is obtained by starting from space--time duality of a Minkowski space,
and a simple Lorentz transformation, that are basic to the theory of special relativity.
The starting point is an adapted version of the continuity equation, namely
\begin{equation*}\label{mod_continuity}
	\partial_t u + \dive_{(t,\mathbf{x})}\, \mathbf{w} = f\qquad(\tau>0),
\end{equation*}
where $\dive_{(t,\mathbf{x})}$ is the scalar product of the operator $(i\sqrt{\tau}\,\partial_t,\partial_{x_1},\dots,\partial_{x_n})$ 
against the extended $(n+1)-$dimensional flux $\mathbf{w}$.
Assuming the extended Fourier's relation
\begin{equation*}
	\mathbf{w}=-a\,\grad_{(t,\mathbf{x})} u,
\end{equation*}
where $\grad_{(t,\mathbf{x})}=(i\sqrt{\tau}\,\partial_t,\grad_{\mathbf{x}})$, we infer
\begin{equation*}\label{scalar_rhd}
	\tau \partial_{t}\!\left(a\,\partial_t u\right) + \partial_t u
		= \dive_{\mathbf{x}}\!\left(a\,\grad_{\mathbf{x}} u\right)+f,
\end{equation*}
which coincides with \eqref{damp_rd} when $a\equiv 1$.
However, the latter equation give rise to significant conceptual issues that makes the theory somewhat controversial.
Among others, some quantities into play are described by complex numbers, with values involving imaginary
``densities'', which are hard to be interpreted.

Finally, let us determine an intermediate form somewhat in between \eqref{damp_rd} and \eqref{scalar_hhd}. 
Let us denote by $\tau_1$ and $\tau_2$ the parameters for \eqref{delayed_continuity} and \eqref{mc_scalar}, respectively.
Combining the delayed version of the continuity equation and the Maxwell--Cattaneo's law
\begin{equation*}
	\left\{\begin{aligned}
	&\tau_1 \partial_{tt} u +\partial_t u + \dive_{\mathbf{x}}\, \mathbf{v} = f,\\
	&\tau_2 \partial_t \mathbf{v} +\mathbf{v} + a\,\grad_{\mathbf{x}} u = 0.
	\end{aligned}\right.
\end{equation*}
Differentiating the first equation with respect to $t$, taking the divergence with respect to $\mathbf{x}$ of the second equation and 
subtracting, we obtain the one-field equation for $u$
\begin{equation*}
	\tau_1\tau_2 \partial_{ttt} u +(\tau_1+\tau_2) \partial_{tt} u +\partial_t\left(u-\tau_2 f\right)
		= \dive_{\mathbf{x}}\left(a\,\grad_{\mathbf{x}} u\right) + f .
\end{equation*}
In the regime of product $\tau_1\tau_2$ small with respect to the other 0-th/1-st order terms in $\tau_1$ and $\tau_2$,
the third order time derivative can be disregarded (if bounded), thus giving raise to the hyperbolic equation
\begin{equation}\label{twodelays}
	\tau \partial_{tt} u +\partial_t\left(u-\sigma f\right)= \dive_{\mathbf{x}}\left(a\,\grad_{\mathbf{x}} u\right) + f .
\end{equation}
where $\tau:=\tau_1+\tau_2$ and $\sigma:=\tau_2$.
In particular, note that $0\leq \sigma\leq \tau$ for any choice of non-negative $\tau_1$ and $\tau_2$.

\subsection{Reaction-diffusion by means of PDE systems}
\label{subsec:systems}

Passing to vector-valued density function $\mathbf{u}\in\mathbb{R}^p$, some modifications have to be taken into account.
First of all, the vectorial form of the continuity equation becomes
\begin{equation}\label{vector_continuity}
	\partial_t \mathbf{u} + \Dive_{\mathbf{x}}\, \mathbf{V}=\mathbf{f} ,
\end{equation}
where $\Dive$ denotes the divergence operator applied to each row of the matrix $\mathbf{V}$,
and $\mathbf{f}$ is some given vector-valued function.

Again, some additional relations coupling the dynamical variables $u$ and $\mathbf{V}$
are required to close the system.
As before, these could be of different nature.
Denoting by $\Grad_{\mathbf{x}}$ the jacobian operator and having in mind the Fourier's law,
we can conceive a relation of the following form
\begin{equation*}\label{vector_genfourier}
	\begin{aligned}
	\mathbf{V} & = \textrm{linear functional applied to } \,\Grad_{\mathbf{x}} \mathbf{u} \\
			& = - \mathbb{A}\,\Grad_{\mathbf{x}} \mathbf{u}.
	\end{aligned}
\end{equation*}
for some ($4^{\textrm{th}}$-order) tensor-valued function $\mathbb{A}$.
Coupling with \eqref{vector_continuity}, the above identity gives the {\it (parabolic) reaction-diffusion system}
\begin{equation}\label{vector_par_rd}
	\partial_t \mathbf{u} = \Dive_{\mathbf{x}}\left(\mathbb{A}\,\Grad_{\mathbf{x}} \mathbf{u}\right)+\mathbf{f},
\end{equation}
which can be regarded as the vectorial extension of the scalar equation \eqref{scalar_par_rd}.

As in Subsection \ref{subsec:scalar}, we may search for alternatives to the Fourier's law,
the first being the Maxwell--Cattaneo's law.
In vectorial version, this reads as
\begin{equation*}\label{vector_maxwellcattaneo}
	\tau\partial_t \mathbf{V} +\mathbf{V} = -\mathbb{A}\,\Grad_{\mathbf{x}} \mathbf{u}.
\end{equation*}
Of course, the latter equality can be generalized to the (more realistic) case in which any line
of the flux matrix $\mathbf{V}$ has a different delay $\tau_1,\dots,\tau_p$.
However, for the sake of simplicity, we will mainly concentrate on the case of a single time-scale $\tau$. 

Coupling with the continuity equation \eqref{vector_continuity}, we end up with the {\it (hyperbolic) reaction-diffusion system with relaxation}
\begin{equation*}\label{vector_hyprd}
	\left\{\begin{aligned}
	\partial_t \mathbf{u} + \Dive_{\mathbf{x}}\mathbf{V} & = \mathbf{f},\\
	\tau\partial_t \mathbf{V} + \mathbb{A}\,\Grad_{\mathbf{x}} \mathbf{u} & = - \mathbf{V}.
	\end{aligned}\right.
\end{equation*}
Applying $\partial_t$ to the first equation, $\Dive_{\mathbf{x}}$ to the second
and taking the difference, we deduce the {\it one-field system}
\begin{equation}\label{vector_hyp_rd}
	\tau\partial_{tt} \mathbf{u} +\partial_t\left(\mathbf{u}-\tau\mathbf{f}\right)
	=\Dive_{\mathbf{x}}\left(\mathbb{A}\,\Grad_{\mathbf{x}} \mathbf{u}\right) + \mathbf{f}.
\end{equation}
The hyperbolic system \eqref{vector_hyp_rd} can be understood as a possible singular perturbation 
of the parabolic limit system \eqref{vector_par_rd}.

Alternatively, we can follow the strategy previously proposed considering a delayed continuity equality,
which ends up in the {\it (hyperbolic) reaction-diffusion system with damping}
\begin{equation}\label{vector_damp_rd}
	\tau\partial_{tt} \mathbf{u}+\partial_t \mathbf{u}
		= \Dive_{\mathbf{x}}\bigl(\mathbb{A}\,\Grad_{\mathbf{x}} \mathbf{u}\bigr)+\mathbf{f}.
\end{equation}
to be regarded as the vectorial version of \eqref{damp_rd}.

In order to derive a sort of interpolation between \eqref{vector_hyp_rd} and  \eqref{vector_damp_rd}, 
we follow the strategy proposed in deducing equation \eqref{twodelays}, that is considering delays
in both continuity identity and flux constitutive equality, with small relaxation times $\tau_1$ and $\tau_2$, so that
the term with the product $\tau_1\tau_2$ can be formally disregarded.
In addition, restricting the attention to
\begin{equation*}
	\mathbb{A}=\textrm{constant} \qquad\textrm{and}\qquad \mathbf{f}=\mathbf{f}(\mathbf{u}),
\end{equation*}
we end up with the system
\begin{equation}\label{finalform1}
	\tau\partial_{tt}\mathbf{u}+\partial_t \left\{\mathbf{u}-\sigma \mathbf{f}(\mathbf{u})\right\}
		= \Dive_{\mathbf{x}}\left\{\mathbb{A}\,\Grad_{\mathbf{x}} \mathbf{u}\right\}+\mathbf{f}(\mathbf{u}),
\end{equation}
Later on, it will be transparent how the apparently harmless term $\sigma\textrm{d}\mathbf{f}(\mathbf{u})$, negligible for $\sigma$ small,
may affect the transient dynamics and plays a crucial role also in the long run.

In the class described by system \eqref{finalform1}, there are some significant limiting regimes, with respect to 
the values of the parameters $\tau$ and $\sigma\in[0,\tau]$:
\begin{itemize}
\item[\bf i.] $\sigma=\tau=0$ (undelayed continuity/undelayed flux): 
\begin{equation*}
	\partial_t \mathbf{u}= \Dive_{\mathbf{x}}\left\{\mathbb{A}\,\Grad_{\mathbf{x}} \mathbf{u}\right\}+\mathbf{f}(\mathbf{u});
\end{equation*}
\item[\bf ii.]  $\sigma=0$, $\tau>0$ (delayed continuity/undelayed flux): 
\begin{equation*}
	\tau\partial_{tt}\mathbf{u}+\partial_t \mathbf{u}
		= \Dive_{\mathbf{x}}\left\{\mathbb{A}\,\Grad_{\mathbf{x}} \mathbf{u}\right\}+\mathbf{f}(\mathbf{u});
\end{equation*}
\item[\bf iii.]  $\sigma=\tau>0$ (undelayed continuity/delayed flux): 
\begin{equation*}
	\tau\partial_{tt}\mathbf{u}+\partial_t\!\left\{\mathbf{u}-\tau \mathbf{f}(\mathbf{u})\right\}
		= \Dive_{\mathbf{x}}\left\{\mathbb{A}\,\Grad_{\mathbf{x}} \mathbf{u}\right\}+\mathbf{f}(\mathbf{u}).
\end{equation*}
\end{itemize}

Additional specifications can be required on the zero-th order term $\mathbf{f}$ to add structure to the whole system.
In the scalar case, any continuous function $f$ has a smooth primitive, producing a corresponding potential $W$, i.e. $W'=-f$.
Differently, when the dimension is strictly greater than 1, additional constraints are needed
in order to make this requirement to be satisfied.
Specifically, for smooth functions, a necessary condition for the existence of a {\it potential function}
$W\,:\,\mathbb{R}^p\to\mathbb{R}^p$ such that
\begin{equation}\label{potential}
	\grad_{\mathbf{u}} W (\mathbf{u}) = -\mathbf{f}(\mathbf{u}),
\end{equation}
is requiring that the jacobian matrix $\textrm{d}\mathbf{f}$ of $\mathbf{f}$ is symmetric, that is
\begin{equation}\label{symm_jac}
	\textrm{d}\mathbf{f}(\mathbf{u})^{\top}=\textrm{d}\mathbf{f}(\mathbf{u}).
\end{equation}
Such condition is also sufficient if the domain for the variable $\mathbf{u}$ is simply connected or star-shaped. 

Assuming the symmetry condition \eqref{symm_jac}, system \eqref{finalform1} is endowed with a natural {\it Lyapunov functional},
i.e. a global function which is not-increasing along any given trajectory $t\mapsto \mathbf{u}(\cdot,t)$ .
To simplify the formalism, we concentrate on the one-dimensional spatial case, limiting ourselves to
\begin{equation}\label{hyperpar}
	\tau\partial_{tt}\mathbf{u}+\partial_t \left\{\mathbf{u}-\sigma \mathbf{f}(\mathbf{u})\right\}
		= \mathbf{A}\,\partial_{xx} \mathbf{u} +\mathbf{f}(\mathbf{u}),
\end{equation}
For $\tau=\sigma=0$, we obtain the standard parabolic reaction-diffusion system
\begin{equation}\label{parpar}
	\partial_t \mathbf{u} = \mathbf{A}\,\partial_{xx} \mathbf{u} + \mathbf{f}(\mathbf{u})
\end{equation}
Property \eqref{potential} guarantees the presence of a {\it variational structure}: the functional
\begin{equation*}
	\mathcal{E}_0[\mathbf{u}]:=\int_{\mathbb{R}} \Bigl\{\tfrac12\mathbf{A}\partial_x \mathbf{u}\cdot \partial_x \mathbf{u}+W(\mathbf{u})\Bigr\}dx,
\end{equation*}
together with some appropriate integrability conditions at $\pm\infty$, is a Lyapunov functional for the system \eqref{parpar}.
Indeed, multiplying by $\partial_t\mathbf{u}$ and integrating by parts, there holds
\begin{equation*}
	\frac{d}{dt}\mathcal{E}_0[\mathbf{u}]+\int_{\mathbb{R}} |\partial_t \mathbf{u}|^2\,dx=0,
\end{equation*}
exhibiting a {\it dissipative property} for $\mathcal{E}_0$, playing the role of an energy functional.

Similar considerations can be done also in the case \eqref{hyperpar}, giving raise
to a differential equality for the {\it modified energy}
\begin{equation*}
	\mathcal{E}_{\tau}[\mathbf{u}]:=\tfrac12\tau|\partial_t\mathbf{u}|^2+\mathcal{E}_0[\mathbf{u}].
\end{equation*}
Then, setting $\mathbf{Q}_\sigma:=\mathbf{I}-\sigma\textrm{d}\mathbf{f}(\mathbf{u})$, there holds
\begin{equation*}
	\frac{d}{dt}\mathcal{E}_{\tau}[\mathbf{u}]+\int_{\mathbb{R}} \mathbf{Q}_\sigma\partial_t \mathbf{u}\cdot \partial_t\mathbf{u}\,dx=0.
\end{equation*}
Again, choosing $\sigma\geq 0$ sufficiently small so that $\mathbf{Q}_\sigma>0$, dissipation is transparent.

\subsection{Three basic numerical schemes in one space dimension}
\label{subsec:algorithms}

For $\mathbf{u}\in\mathbb{R}^p$ and in one space dimension, the tensor $\mathbb{A}$ reduces now to a $p\times p$ matrix $\mathbf{A}$,
that is $\mathbb{A}=\mathbf{A}=(a_{1\ell }^{i1})$ since two of the four indeces are now fixed and equal to $1$.
For the sake of simplicity, we limit ourselves to the case $\mathbf{A}=a\,\mathbf{I}$ for some constant $a>0$.
Hence, we consider the system in one space-dimension
\begin{equation}\label{parhyp_rd_1d}
	\tau\partial_{tt} \mathbf{u} +\partial_t \left\{\mathbf{u}-\sigma \mathbf{f}(\mathbf{u})\right\}
		= a\partial_{xx} \mathbf{u}+\mathbf{f}(\mathbf{u}).
\end{equation}
Let us stress once more that the idea is not to consider hyperbolic models as perturbations of the limiting parabolic ones,
but rather to explore numerical approximation of the hyperbolic equations regarded as intriguing models on their own with 
different properties, with particular care to the transient behavior.
Later on, we will test and compare the numerical schemes with specific attention to their capability of providing
precise approximations of the propagation speed of the special solutions called {\it fronts}.

\subparagraph{First-order reduction algorithm}
System \eqref{parhyp_rd_1d} has an immediate numerical description, obtained by rewriting it in first-order form as
\begin{equation}\label{1storderred}
	\left\{\begin{aligned}
	\partial_t \mathbf{u} & = \mathbf{v},\\
	\tau\,\partial_{t} \mathbf{v} & = a\,\partial_{xx}\mathbf{u}+\mathbf{f}(\mathbf{u})-\left\{\mathbf{I}-\sigma\,\textrm{d}\mathbf{f}(\mathbf{u})\right\}\mathbf{v}.
	\end{aligned}\right.
\end{equation}
Firstly, we discretize the spatial part by introducing a uniform mesh with step $\mathrm{dx}$, 
\begin{equation}\label{1stord_semidiscrete} 
	\left\{\begin{aligned}
	\frac{d\mathbf{u}_j}{dt}&=\mathbf{v}_j\\
	\tau\,\frac{d\mathbf{v}_j}{dt}&= \frac{a}{\mathrm{dx}^2}\,\left(\mathbf{u}_{j+1}-2\mathbf{u}_{j}+\mathbf{u}_{j-1}\right)
		+\mathbf{f}(\mathbf{u}_j)-\left\{\mathbf{I}-\sigma\,\textrm{d}\mathbf{f}(\mathbf{u}_j)\right\}\mathbf{v}_j
	\end{aligned}\right.
\end{equation}
Then, a subsequent time-discretization, that can be performed in different ways, is applied.
To start with, we choose an implicit-explicit scheme (IMEX), limiting the implicit description to the linear part of the system, so that 
\begin{equation*}
	\left\{\begin{aligned}
	\frac{\mathbf{u}_j^{n+1}-\mathbf{u}_j^{n}}{\mathrm{dt}}&=\mathbf{v}_j^{n+1}\\
	\tau\,\frac{\mathbf{v}_j^{n+1}-\mathbf{v}_j^{n}}{\mathrm{dt}}&= \frac{a}{\mathrm{dx}^2}\,\left(\mathbf{u}_{j+1}^{n+1}-2\mathbf{u}_{j}^{n+1}+\mathbf{u}_{j-1}^{n+1}\right)
		+\mathbf{f}(\mathbf{u}_j^{n})-\mathbf{v}_j^{n+1}+ \sigma \textrm{d}\mathbf{f}(\mathbf{u}_j^{n})\mathbf{v}_j^{n}
	\end{aligned}\right.
\end{equation*}
which gives the {\it first-order (reduction) algorithm}
\begin{equation}\label{imex1}
	\left\{\begin{aligned}
	\mathbf{u}_j^{n+1}-\mathrm{dt}\,\mathbf{v}_j^{n+1}&=\mathbf{u}_j^{n}\\
	\alpha\left(-\mathbf{u}_{j+1}^{n+1}+2\mathbf{u}_{j}^{n+1}-\mathbf{u}_{j-1}^{n+1}\right)+\left(1+ \beta\right) \mathbf{v}_j^{n+1}
	& = \mathbf{v}_j^{n} + \beta\,\mathbf{f}(\mathbf{u}_j^{n}) \\
	&\quad + \sigma\beta\,\textrm{d}\mathbf{f}(\mathbf{u}_j^{n})\mathbf{v}_j^{n},
	\end{aligned}\right.
\end{equation}
where $\alpha:=a\,\mathrm{dt}/\tau\mathrm{dx}^2$, $\beta:=\mathrm{dt}/{\tau}$.
Solving such an implicit-explicit algorithm furnishes the numerical approximation of the real solution
\begin{equation*}\label{imex2}
	\begin{pmatrix} \mathbf{u}^{n+1} \\ \mathbf{v}^{n+1} \end{pmatrix}
	=\mathbf{A}^{-1} \begin{pmatrix} \mathbf{u}^{n} \\ \mathbf{v}^{n} + \beta\,\mathbf{f}(\mathbf{u}^{n})
		+ \sigma\beta\,\textrm{d}\mathbf{f}(\mathbf{u}^{n})\mathbf{v}^{n}\,\mathrm{dt} \end{pmatrix}
\end{equation*}
where $\mathbf{A}$ describes the coefficients of the left-hand side matrix in \eqref{imex1}.

\subparagraph{Li\'enard-type algorithm}
A second type of algorithm is inspired by the so-called {\it Li\'enard second order equation} which is
\begin{equation*}
	\tau\frac{d^2 u}{dt^2}+g(u)\frac{du}{dt}+h(u)=0.
\end{equation*}
for some given functions $g$ and $h$.
The above equation can be rewritten as a first order system by setting
\begin{equation*}
	\tau\frac{du}{dt}=v-G(u),\qquad \frac{dv}{dt}=-h(u)
\end{equation*}
where $G$ is a primitive of the function $g$.
Applied to system \eqref{parhyp_rd_1d}, let us consider an algorithm,
which will be later named {\it Li\'enard-type algorithm},
based on the decomposition
\begin{equation}\label{lienard}
	\left\{\begin{aligned}
	\tau \partial_{t} \mathbf{u} & = \mathbf{v}-\mathbf{u}+ \sigma \mathbf{f} (\mathbf{u}),\\
	\partial_t \mathbf{v} & = a\,\partial_{xx} \mathbf{u}+\mathbf{f}(\mathbf{u}).
	\end{aligned}\right.
\end{equation}
As before, discretizing with respect to a mesh with step $\mathrm{dx}$, we infer
\begin{equation}\label{lienard_semidiscrete}
	\left\{\begin{aligned}
	\tau\,\frac{d\mathbf{u}_j}{dt}&=\sigma \mathbf{f} (\mathbf{u}_j) -\mathbf{u}_j + \mathbf{v}_j\\
	\frac{d\mathbf{v}_j}{dt}&= \frac{a}{\mathrm{dx}^2}\,\left(\mathbf{u}_{j+1}-2\mathbf{u}_{j}+\mathbf{u}_{j-1}\right)+\mathbf{f}(\mathbf{u}_j).
	\end{aligned}\right.
\end{equation}
At the continuous level, systems \eqref{1storderred} and \eqref{lienard}, and the corresponding semi-discrete
algorithms, viz. systems \eqref{1stord_semidiscrete}  and \eqref{lienard_semidiscrete}, are completely equivalent, the difference
being only in the choice of the variable $\mathbf{v}$. 

Distinctions emerge in the subsequent step, where the time discretiza\-ti\-on is taken into account
and the difference between linear (implicit) vs nonlinear (explicit) discretizations emerges.
On top of that, we observe that the Li\'enard-type algorithm does not require an explicit computation
of the jacobian matrix $\textrm{d}\,\mathbf{f}$ at the value $\mathbf{u}^n_j$;
hence, in principle, it could be considered also for less smooth reaction term $\mathbf{f}$.

Proceeding in the same spirit as above, we infer
\begin{equation*}
	\left\{\begin{aligned}
	\tau\frac{\mathbf{u}_j^{n+1}-\mathbf{u}_j^{n}}{\mathrm{dt}}&=\sigma \mathbf{f} (\mathbf{u}_j^n) -\mathbf{u}_j^{n+1} + \mathbf{v}_j^{n+1}\\
	\frac{\mathbf{v}_j^{n+1}-\mathbf{v}_j^{n}}{\mathrm{dt}}&= \frac{a}{\mathrm{dx}^2}\,\left(\mathbf{u}_{j+1}^{n+1}-2\mathbf{u}_{j}^{n+1}+\mathbf{u}_{j-1}^{n+1}\right)
		+\mathbf{f}(\mathbf{u}_j^{n})
	\end{aligned}\right.
\end{equation*}
from which we obtain the IMEX linear system
\begin{equation}\label{imex1_lienard}
	\left\{\begin{aligned}
	(1+\beta)\mathbf{u}_j^{n+1}-\beta \mathbf{v}_j^{n+1}
		&=\mathbf{u}_j^{n}+\beta \sigma \mathbf{f} (\mathbf{u}_j^n) \\
	\alpha \left(-\mathbf{u}_{j+1}^{n+1}+2\mathbf{u}_{j}^{n+1}-\mathbf{u}_{j-1}^{n+1}\right)+\mathbf{v}_j^{n+1}
		&=\mathbf{v}_j^{n}+\mathbf{f}(\mathbf{u}_j^{n})\mathrm{dt}
	\end{aligned}\right.
\end{equation}
with $\alpha:=a\mathrm{dt}/\mathrm{\mathrm{dx}^2}$ and $\beta:=\mathrm{dt}/{\tau}$.
The solution of such an iteration provides the numerical approximation of the solution $(\mathbf{u},\mathbf{v})$
\begin{equation*}\label{imex2_lienard}
	\begin{pmatrix} \mathbf{u}^{n+1} \\ \mathbf{v}^{n+1} \end{pmatrix}
	=\mathbf{A}^{-1} \begin{pmatrix} \mathbf{u}_j^{n}+\beta \sigma \mathbf{f} (\mathbf{u}_j^n) \\
								\mathbf{v}_j^{n}+\mathbf{f}(\mathbf{u}_j^{n})\mathrm{dt} \end{pmatrix}
\end{equation*}
where $\mathbf{A}$ describes the coefficients of the left-hand side matrix in \eqref{imex1_lienard}.

\subparagraph{Kinetic algorithm}

A third viable algorithm is limited to the special case $\sigma=\tau>0$.
In such a situation,  let us start back from the derivation of the model, 
i.e. the coupling of the balance law together with the Maxwell--Cattaneo's relation,
\begin{equation*}
	\partial_t \mathbf{u} + \partial_x \mathbf{v} = \mathbf{f} (\mathbf{u}),\qquad
	\tau\,\partial_t \mathbf{v} + \mathbf{v} = - a\,\partial_x \mathbf{u}.
\end{equation*}
Here, $\tau$ and $a$ can be considered as diagonal matrices with elements $(\tau_1,\dots,\tau_n)$ and $(a_1,\dots,a_n)$,
 with components $\tau_i$ and $a_i$ which are considered possibly different one from the other.
Therefore, we end up with the system
\begin{equation}\label{relax_1d_compo}
	\left\{\begin{aligned}
	\partial_t u_i + \partial_x v_i & = f_i (u_1,\dots,u_n),\\
	\tau_i\, \partial_t v_i + a_i\,\partial_x u_i & = - v_i.
	\end{aligned}\right.
\end{equation}
The coupling is due to the presence of the term $\mathbf{f}=(f_1,\dots,f_n)$ in the first equation.

The coefficients of the principal part of the differential operator  at the left-hand side of \eqref{relax_1d_compo} are described
by the block-diagonal matrix $\mathbf{A}=\textrm{blockdiag}(\mathbf{A}_1,\dots,\mathbf{A}_n)$ with
\begin{equation*}
	\mathbf{A}_i:=\begin{pmatrix} 0 & 1 \\ a_i/\tau_i &0 \end{pmatrix}
	\qquad\qquad i=1,\dots,n.
\end{equation*}
 Therefore, the eigenvalues of the matrix $\mathbf{A}$ with size $2n$,  are given by the roots of the polynomial
\begin{equation*}
	p(\lambda)=\det(\mathbf{A}-\lambda\mathbf{I})=\prod_{i=1}^{n}\left(\lambda^2-\varrho_i^2\right),
\end{equation*}
where $\varrho_i:=\sqrt{{a_i}/{\tau_i}}$, are $\lambda=\pm\varrho_i$ for $i=1,\dots, n$.

Introducing the diagonal variables $(\mathbf{r},\mathbf{s})=(r_1,\dots,r_n,s_1,\dots,s_n)$, defined by
\begin{equation*}
	r_i:=\frac{1}{2}\left(u_i-\frac{v_i}{\rho_i}\right),\qquad
	s_i:=\frac{1}{2}\left(u_i+\frac{v_i}{\varrho_i}\right),
\end{equation*}
system \eqref{relax_1d_compo} becomes
\begin{equation*}
	\left\{\begin{aligned}
	\partial_{t} r_i - \varrho_i \partial_{x} r_i
		&= \frac{1}{2\tau_i}(-r_i+s_i)+\frac12 f_i(r_1+s_1,\dots,r_n+s_n),\\
	\partial_{t} s_i + \varrho_i \partial_{x} s_i
		&= \frac{1}{2\tau_i}(+r_i-s_i)+\frac12 f_i(r_1+s_1,\dots,r_n+s_n).
	\end{aligned}\right.
\end{equation*}
As before, we firstly consider a spatial discretization with a uniform mesh  of step $\mathrm{dx}$.
 Taking into account the up-wind nature of the model, we obtain
\begin{equation*}
	\left\{\begin{aligned}
	\frac{d r_{i,j}}{dt} - \varrho_i \frac{r_{i,j+1}-r_{i,j}}{\mathrm{dx}}
		+ \frac{1}{2\tau_i}(+r_{i,j}-s_{i,j}) &= \frac12 f_i(r_1+s_1,\dots,r_n+s_n),\\
	\frac{d s_{i,j}}{dt} + \varrho_i \frac{s_{i,j}-s_{i,j-1}}{\mathrm{dx}}
		+ \frac{1}{2\tau_i}(-r_{i,j}+s_{i,j}) &=\frac12 f_i(r_1+s_1,\dots,r_n+s_n).
	\end{aligned}\right.
\end{equation*}
Next, we follow the same strategy of the IMEX algorithm, that is we discretize implicitly only the linear
part of the system.
Thus, we infer
\begin{equation*}
	\left\{\begin{aligned}
	\frac{r_{i,j}^{n+1}-r_{i,j}^{n}}{\mathrm{dt}}-\frac{\varrho}{\mathrm{dx}}\bigl(r_{i,j+1}^{n+1}-r_{i,j}^{n+1}\bigr)
		&- \frac1{2\tau_i}\bigl(-r_{i,j}^{n+1}+s_{i,j}^{n+1}\bigr)=\frac12 f_i(\mathbf{r}^n+\mathbf{s}^n),\\
	\frac{s_{i,j}^{n+1}-s_{i,j}^{n}}{\mathrm{dt}}+\frac{\varrho}{\mathrm{dx}}\bigl(s_{i,j}^{n+1}-s_{i,j-1}^{n+1}\bigr)
		&- \frac1{2\tau_i}\bigl(r_{i,j}^{n+1}-s_{i,j}^{n+1}\bigr)=\frac12 f_i(\mathbf{r}^n+\mathbf{s}^n),
	\end{aligned}\right.
\end{equation*}
that gives
\begin{equation*}\label{KinAlgoSystem1}
	\left\{\begin{aligned}
	\left(1+\alpha_i+\beta_i\right)r_{i,j}^{n+1}-\alpha_i r_{i,j+1}^{n+1} - \beta_i s_{i,j}^{n+1}&=r_{i,j}^{n}+\frac12 f_i(\mathbf{r}^n+\mathbf{s}^n)\,\mathrm{dt},\\
	- \beta_i r_{i,j}^{n+1}-\alpha_i s_{i,j-1}^{n+1}+ \left(1+\alpha_i + \beta_i \right)s_{i,j}^{n+1}&=s_{i,j}^{n}+\frac12 f_i(\mathbf{r}^n+\mathbf{s}^n)\,\mathrm{dt},
	\end{aligned}\right.
\end{equation*}
where $\alpha_i=\varrho_i\,\textrm{dt}/\textrm{dx}$, $\beta_i=\textrm{dt}/2\tau_i$.
Again, denoting by $\mathbf{A}$ the coefficients' matrix of the couple $(\mathbf{r},\mathbf{s})$ in the above system,
we obtain the iteration formula
\begin{equation*}\label{KinAlgoSystem2}
	\begin{pmatrix} \mathbf{r}^{n+1}\\ \mathbf{s}^{n+1} \end{pmatrix} = 
	\mathbf{A}^{-1} \begin{pmatrix}
		\mathbf{r}^{n} + \mathbf{f}(\mathbf{r}^n+\mathbf{s}^n)\,\mathrm{dt}/2 \\
		\mathbf{s}^{n} + \mathbf{f}(\mathbf{r}^n+\mathbf{s}^n)\,\mathrm{dt}/2 \end{pmatrix}\,,
\end{equation*}
 defining the mapping at the base of the numerical algorithm.

\section{Some waves are better than others}
\label{sect:somewaves}

\setcounter{equation}{0}
\renewcommand{\theequation}{2.\arabic{equation}}

 In this Section the attention moves towards a class of particularly significant special solutions:
the {\it traveling waves}.
 Such solutions are indeed supported by hyperbolic reaction-diffusion system  corresponding to scalar
parabolic reaction-diffusion equations for both monostable and bistable reaction terms.
 Moreover, we focus on a special class of waves, called {\it propagation fronts} explored in details in the case 
of a bistable reaction term.
Special cases where the speed of propagation can be explicitly computed are also provided.
 A detailed discussion on the monostable case can be found in \cite{Hade88} (see also \cite{BouiCalvNadi14}).

\subsection{Traveling waves}
\label{subsect:twaves}

Among the infinitely many solutions of a partial differential equations, some solutions exhibits  usually an augmented ``stability'',
 inherited by the additional amount of internal symmetry.
A recurrent type of such kind of solutions are the so-called {\it   planar} traveling waves  (or simply 
{\it traveling waves}), i.e. solutions of the form
\begin{equation}\label{planarfront}
	\mathbf{u}(\mathbf{x},t):=\boldsymbol{\phi}(\mathbf{k}\cdot\mathbf{x}-c\, t)
\end{equation}
for some unitary vector $\mathbf{k}$.
Here $\boldsymbol{\phi}$ is called the {\it profile of the wave} and $c$ its {\it propagation speed}. 

For such  special solutions, PDEs are reduced to ODEs  with unknowns depending on the scalar variable
$\xi:=\mathbf{k}\cdot\mathbf{x}-c\, t$  and for a value $c$ to be determined together with the function $\boldsymbol{\phi}$.
As an example, inserting the ansatz \eqref{planarfront} in \eqref{finalform1} and noticing that
\begin{equation*}
	\begin{aligned}
	\Grad_{\mathbf{x}} \mathbf{u}&=\frac{d\boldsymbol{\phi}}{d\xi}\otimes \mathbf{k},\\
	\Dive_{\mathbf{x}}\left\{\mathbb{A}\,\left(\frac{d\boldsymbol{\phi}}{d\xi}\otimes \mathbf{k}\right)\right\}
		&=\mathbb{A}\,\Dive_{\mathbf{x}}\left(\frac{d\boldsymbol{\phi}}{d\xi}\otimes \mathbf{k}\right)
		=\mathbb{A}\,\left(\frac{d^2 \boldsymbol{\phi}}{d\xi^2}\otimes \mathbf{k}\right),
	\end{aligned}
\end{equation*}
we end up with an ODE for the profile $\phi$, parametrized by the velocity $c$,
\begin{equation*}
	 \mathbb{A}\,\left(\frac{d^2 \boldsymbol{\phi}}{d\xi^2}\otimes \mathbf{k}\right)+c^2 \tau\,\frac{d^2 \boldsymbol{\phi}}{d\xi^2}
	 	+c \left\{\mathbf{I}-\sigma d\mathbf{f}(\boldsymbol{\phi})\right\}\frac{d\boldsymbol{\phi}}{d\xi}+\mathbf{f}(\boldsymbol{\phi})=0.
\end{equation*}
Since the above system is autonomous, the profile is determined up to translations.
 In particular, translation $\boldsymbol{\phi}_\delta:=\boldsymbol{\phi}(\cdot-\delta)$ with $\delta\in\mathbb{R}$
of a given traveling wave $\boldsymbol{\phi}=\boldsymbol{\phi}(\cdot)$ is itself a traveling wave solution for the same  equation.
Such properties have an immediate consequence: the derivative of $\boldsymbol{\phi}$ with respect to its argument
is an eigenfunction for the corresponding linearized operator at $\boldsymbol{\phi}$ relative to the eigenvalue $\lambda=0$.
This influences the stability properties of the wave, dictating the fact that, at most, {\it orbital stability} could be expected,
meaning convergence of small perturbations to the manifold $ \Phi:=\{\boldsymbol{\phi}_\delta\,:\,\delta\in\mathbb{R}\}$.
Presence/absence of an {\it asymptotic phase} --viz. convergence to a definite element of the manifold  $\Phi$--
is the (natural) subsequent issue.

Depending on specific properties of the profile function $\boldsymbol{\phi}$, different names are associated to traveling waves:
\begin{itemize}
\item[\bf i.] 	if $\boldsymbol{\phi}$ converges to some asymptotic states $\boldsymbol{\phi}_\pm$ (which are necessarily
two equilibria of the model) with $\boldsymbol{\phi}_-\neq \boldsymbol{\phi}_+$, the solution is called a {\it front};
\item[\bf ii.]	if $\boldsymbol{\phi}$ converges to the same asymptotic state $\overline{\boldsymbol{\phi}}$
( again, equilibrium of the model), the solution is said to be a {\it pulse};
\item[\bf iii.]	if $\boldsymbol{\phi}$ is periodic, the solution is a {\it wave-train}.
\end{itemize}

In the state space, the three configurations correspond, respectively, to the presence of a heteroclinic orbit,
a homoclinic orbit, a cycle.
From now on, we focus on  the analysis of fronts; also, we restrict the attention to the spatial one-dimensional case.
A further reduction concerns with the size of the vector $\mathbf{u}$ which is, from the time being,
regarded as a scalar quantity $u$, thus restricting the attention to the second-order scalar equation
\begin{equation}\label{scalar_damp_rd_1d}
	\tau\,\partial_{tt} u + \partial_t\left\{u-\sigma f(u)\right\} = a\,\partial_{xx} u+f(u).
\end{equation}
where $f$ is an appropriate functions and $\tau, \sigma, a$ are positive constants  with $\sigma\in[0,\tau]$.

\subparagraph{Monostable and bistable nonlinearities}

Following \cite{AronWein78,AronWein75}, we  focus on two types of nonlinearities.
\begin{itemize}
\item[\bf i.] {\sl Monostable.}
	The function $f$ is assumed to be smooth, strictly positive in some fixed interval $(0,1)$,
	negative in $(-\infty,0)\cup(1,+\infty)$, with simple zeros; 
\item[\bf ii.] {\sl Bistable.}
	The function $f$ is assumed to be smooth, strictly positive in some fixed interval
	$(-\infty,0)\cup(\alpha,1)$, negative in $(0,\alpha)\cup(1,+\infty)$, with simple zeros.
\end{itemize}

\begin{figure}[htb]\centering
\includegraphics[width=10cm]{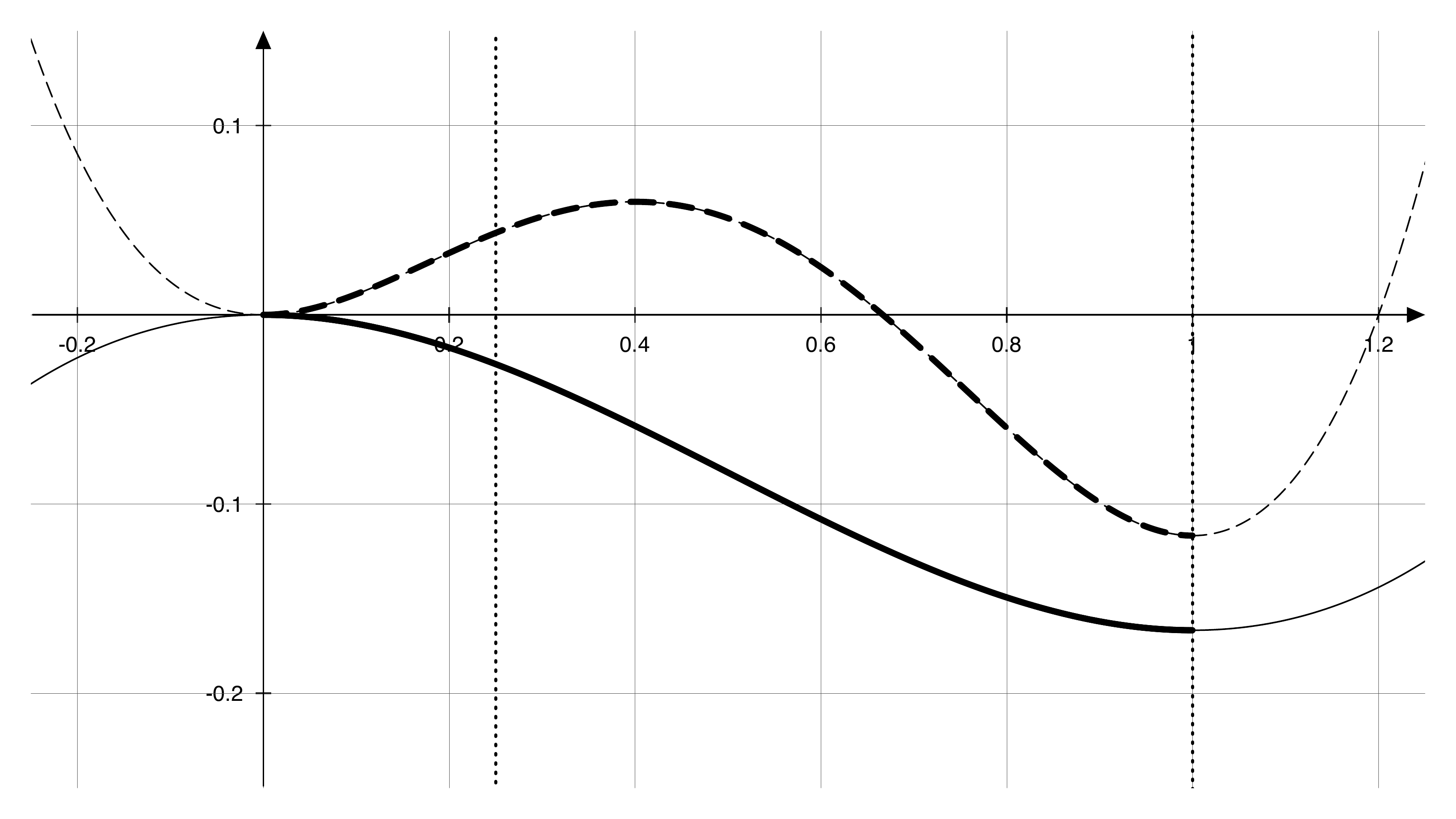}
\caption{The potentials $W$ relative to the functions $f(u)=u(1-u)$ (monostable, continuous)
and $f(u)=\kappa u(1-u)(u-\alpha)$ with $\kappa=7$ and $\alpha=0.4$ (bistable, dashed).}
\label{fig:potentials}
\end{figure}

In both situations, we introduce the corresponding potential 
\begin{equation*}\label{scalarpotential}
	W(u):=-\int_{0}^{u} f(s)\,ds
\end{equation*}
The function $W$ is decreasing for the monostable regime and it has a double-well form for the bistable one
 (see Fig.\ref{fig:potentials}).

The former case, whose prototype is $f(u)\propto u(1-u)$, corresponds to a logistic-type reaction term and it is usually
referred to as {\it Fisher--KPP equation} (using the initials of the names Kolmogorov, Petrovskii and Piscounov).
The potential corresponding to  the logistic function $f(u)=\kappa\,u(1-u)$ is
\begin{equation*}\label{logistic_pot}
	W(u)=\tfrac{1}{6}\kappa\left(-3u^2+2u^3\right),
\end{equation*}
drawn in Figure \ref{fig:potentials} (continuous line).
Different kind of monostable reaction function $f$ are the {\it Gompertz term}, i.e. $f(u)=\kappa\,u \ln u$, and {\it von Bertalanffy term},
i.e. $f(u)=\kappa(u^\mu - u)$ with $\mu\in(0,1)$, corresponding potentials being $W(u)=\kappa u^2(2\ln|u|-1)/4$ and
$W(u)=\kappa\left\{u^2/2-u^{\mu+1}/(\mu+1)\right\}$, respectively.
The main difference is in the location of the tangent line at $u=0$, vertical in the last two cases,
and  playing a crucial role in the statement of existence of propagating fronts.

The latter, whose behaviour is roughly given by the third order polynomial $f(u)\propto u(u-\alpha)(1-u)$ with $\alpha\in(0,1)$,
is called {\it Allen--Cahn equation} (sometimes, also bear the names of {\it Nagumo} and/or {\it Ginzburg--Landau}).
The potential which corresponds to $f(u)=\kappa\,u(u-\alpha)(1-u)$ is
\begin{equation}\label{verhulst_pot}
	W(u)=\tfrac{1}{12}\kappa\,u^2\left\{6\alpha-4(1+\alpha)u+3u^2\right\}.
\end{equation}
The presence of the additional intermediate zero of $f$ given by $\alpha$ emerged in ecological context
where it describes the so-called {\it Allee-type effect}, needed  when cooperation is required for
survival (see \cite{CourBereGasc08} for a detailed description of the topic).

\subsection{Propagating fronts}
\label{subsect:localglobal}

Both monostable and bistable nolinearities share a common crucial feature:
they support existence of heteroclinic traveling waves.

\begin{defn}
A {\it propagating front} is a traveling wave solution for a given PDE system having the special form $u(x,t)=\phi(\xi)$
where $\xi:=x-ct$, connecting two different asymptotic states $\phi(\pm \infty)=\phi_\pm$ with $\phi_-\neq \phi_+$.
\end{defn}

The main goal stems in showing existence of a heteroclinic solution to the corresponding second order differential equation
\begin{equation}\label{twode}
	(a-\tau c^2)\frac{d^2\phi}{d\xi^2} + c\frac{d}{d\xi}\left\{\phi+\sigma \frac{dW}{du}(\phi)\right\} - \frac{dW}{du}(\phi) = 0,
\end{equation}
with boundary conditions $\phi(-\infty)= \phi_-$, $\phi(+\infty)= \phi_+$
 where we assume, for definiteness, $\phi_-=1$ and $\phi_+=0$.

Equivalentlly, the second order differential equation \eqref{twode} can be rewritten as
\begin{equation}\label{twode_system}
	\left\{\begin{aligned}
	\frac{d\phi}{d\xi}&=\psi,\\
	\frac{d\psi}{d\xi}&=\frac{1}{a-\tau c^2}\left\{\frac{dW}{du}(\phi)-c\left[1+\sigma\,\frac{d^2 W}{du^2}(\phi)\right]\psi\right\}
	\end{aligned}\right.
\end{equation}
Next, assume $1+\sigma W''(s)>0$ for any $s$ under consideration, 
 which is indeed satisfied if $\sigma$ is sufficiently small.
Multiplying by $d\phi/d\xi$, we deduce the identity
\begin{equation*}\label{somerelation}
	\frac{d}{d\xi}\left\{\tfrac12(a-\tau c^2)\left(\frac{d\phi}{d\xi}\right)^2-W(\phi)\right\}
	+c\left[1+\sigma\,\frac{d^2 W}{du^2}(\phi)\right]\left(\frac{d\phi}{d\xi}\right)^2=0.
\end{equation*}
Thus, integrating in $\mathbb{R}$, we infer
\begin{equation}\label{speedformula1}
	c=\frac{W(0)-W(1)}{\int_{\mathbb{R}}\left[1+\sigma W''(\phi)\right]\left(\phi'\right)^2 dx}
\end{equation}
From this relation, it is readily observed that the speed $c$ is strictly positive if and only if $W(1)<W(0)$.
In particular, in the monostable case, $\phi_-=0$ is a maximum point and $\phi_+=1$ is a minimum for $W$
and thus $c$ is strictly positive.
Differently, in the bistable case, $W$ is a double-well potential and thus the speed is positive or negative depending
on the depth difference  $W(0)-W(1)$ of the two wells  located at $0$ and $1$.

The starting point in proving existence of propagation fronts is the stability analysis of the singular points of \eqref{twode},
i.e. constant values $\bar u$ with the property $f(\bar u)=0$, with respect to the ordinary differential system
obtained by considering the traveling wave ansatz where the speed $c$ is, for the moment, an external parameter.

Linearizing at $\bar u$ the second order differential equation \eqref{twode}, we infer
\begin{equation}\label{firstorderTw}
	\left\{\begin{aligned}
	\frac{d\phi}{d\xi}&=\psi,\\
	\frac{d\psi}{d\xi}&=\frac{1}{a-\tau c^2}\left\{W''(\bar u)\,\phi-c\left[1+\sigma W''(\bar u)\right]\psi\right\}.
	\end{aligned}\right.
\end{equation}
The  corresponding characteristic polynomial is
\begin{equation*}\label{charpolyTw}
	p(\lambda;\bar u,c):=\frac{1}{(a-\tau c^2)}\left\{(a-\tau c^2)\lambda^2+c\left(1+\sigma W''(\bar u)\right)\lambda-W''(\bar u)\right\}.
\end{equation*}
Thus, setting
\begin{equation}\label{Delta}
	\Delta(\bar u,c):=c^2\left\{1+\sigma W''(\bar u)\right\}^2+4(a-\tau c^2)W''(\bar u)>0,
\end{equation}
the two roots of $p=p(\cdot;\bar u,c)$ are
\begin{equation}\label{proots}
	\lambda_\pm(\bar u,c)=\frac{-c\left\{1+\sigma W''(\bar u)\right\}\pm\sqrt{\Delta(\bar u,c)}}{2(a-\tau c^2)}.
\end{equation}
Since they have opposite signs if $W''(\bar u)>0$, the singular point $(\bar u,0)$ is a saddle point for \eqref{firstorderTw}.
Differently, if $W''(\bar u)<0$ (hence $\bar u$ is unstable with respect to the PDE), the two roots are either complex conjugates
or both real with the same sign, thus they define either a spiral or a node.
Assuming $1+\sigma W''(\bar u)>0$, the spiral and the node are stable (or unstable, respectively) if $c>0$ (or $c<0$, resp.).
Hence, the heteroclinic orbit is a node/saddle connection in the case of Fisher--KPP equation (monostable case) and a saddle/saddle connection
in the case of the Allen--Cahn equation (bistable case) for both the parabolic ($\tau=0$) and the hyperbolic equations ($\tau>0$),
 with relevant consequence in term of the multiplicity of the speeds  $c$.

To fix idea, let us give a closer look to traveling waves with a monotone decreasing profile, that is $\phi_-:=1>0=:\phi_+$.
The opposite case can be deduced by straightforward symmetry arguments.

For the node/saddle connection, the situation is rather complicated.
First of all, we have to restrict the attention to the regimes of the parameter $c$ such that the critical point is an unstable node,
ruling out stable/unstable spirals and stable nodes.
For $W''(\bar u)<0$,  the discriminant $\Delta$,  defined in \eqref{Delta}, distinguishes whether the two roots
of the polynomial $p$ are real or not.
When strictly positive, such roots are real and distinct  and we search for intersection between the two-dimensional
unstable manifold of the critical point $\phi_-=1$ at $-\infty$ and the one-dimensional stable manifold at $\phi_+=0$ at $+\infty$. 
In term of dimensions, the situation is favourable.
 Additional computations show that existence could be provided
for a whole half-line of values for the parameter $c$.
For more details on the monostable case, we refer to  \cite{HadeRoth75} in the parabolic case (i.e. $\sigma=\tau=0$)
and to \cite{BouiCalvNadi14} for the case $\sigma=\tau>0$.

 For the saddle/saddle connection, the one-dimensional manifold of the steady state $\phi_-=1$ has to intersect
the stable manifold of the steady state $\phi_+=0$.
Being the system planar, the corresponding stable and unstable manifolds are one-dimensional and
the intersection of the two manifolds is non-generic,
corresponding to the fact that the speed $c$ has to be appropriately tuned.
This translates into the existence of a specific value of the speed  $c_\ast$ for which the heteroclinic connection emerges.

From now on, we restrict the attention to the bistable case with $W(1)\leq W(0)$ so that $c\geq 0$,  see formula \eqref{speedformula1},
with the exception of some minor deviations from the mainstream dedicated to the monostable case.
In particular, we may restrict the attention to the sub-characteristic  regime, determined by the additional requirement $a-\tau c^2>0$.

Introducing the variable $\zeta:=\xi/\sqrt{a-\tau c^2}$, equation \eqref{twode} becomes simpler, namely
\begin{equation}\label{twode2}
	\frac{d^2\phi}{d\zeta^2} + \gamma\frac{d}{d\zeta}\left\{\phi+\sigma W'(\phi)\right\} - W'(\phi) = 0,
\end{equation}
where
\begin{equation}\label{rescaledspeed}
	c_\tau:=\frac{c}{\sqrt{a-\tau c^2}}.
\end{equation}
Equation \eqref{twode2} can be equivalently rewritten as the first order system 
\begin{equation}\label{twode_1st}
	\left\{\begin{aligned}
	\frac{d\phi}{d\zeta} &=\psi,\\ 
	\frac{d\psi}{d\zeta} &= W'(\phi) - c_\tau\left\{1+\sigma W''(\phi)\right\}\psi
	\end{aligned}\right.
\end{equation}
with asymptotic conditions $(\phi,\psi)(-\infty)=(1,0)$ and $(\phi,\psi)(+\infty)=(0,0)$.
A different first order form for \eqref{twode2} is given by the Li\'enard form
\begin{equation*}\label{twode_lienard}
	\left\{\begin{aligned}
		\frac{d\phi}{d\zeta}&=-c_\tau\left\{\phi+\sigma W'(\phi)\right\}+\chi,\\
		\frac{d\chi}{d\zeta} &= W'(\phi),
	\end{aligned}\right.
\end{equation*}
with asymptotic conditions $(\phi,\chi)(-\infty)=(1,0)$ and $(\phi,\chi)(+\infty)=(0,0)$.

 The simplified form \eqref{twode_1st} for \eqref{twode} is particularly convenient when passing from local to global analysis,
using the {\it rotated vector field property} of system \eqref{twode2}.
The final statement relative to existence of propagating front is reported here, for readers' convenience, as taken from \cite{LattMascPlazSimeXX}.

\begin{thm}\label{theoexists}
Let $W$ be a double-well potential with local minima at $0$ and $1$.
If $\tau>0$, $\sigma\in[0,\tau]$ and $1+\sigma W''(s)>0$ for any $s\in[0,1]$,
then there exists a unique value $c_\ast\in\mathbb{R}$ such that the equation \eqref{twode}
has a monotone increasing solution $\phi$ with asymptotic states $\phi(-\infty)=1$ and $\phi(+\infty)=0$.
\end{thm}

\subsection{Special cases with explicit propagation speeds}
\label{subsect:specialcases}

Next, we focus on three special cases for which an explicit formula is available.
The first one concerns with the case of two wells of equal depth.
Next, we pass to consider the specific case of a third order polynomial reaction term for which explicit formulas
for both the standard parabolic equation and the  damped hyperbolic one can be determined.
Finally, we discuss the case of a piecewise linear reaction function with a jump located at  some intermediate value $\alpha$.

\subparagraph{Two wells of equal depth}

The case of a double-well potential $W$ with wells of equal depth can be treated separately,
since \eqref{speedformula1} indicates that $c_\ast=0$, indipendently on the values of $\sigma\geq 0$.

\begin{prop}
Let $\tau\geq 0$ and $\sigma\in[0,\tau]$.
In addition, let $f=-W'$ with $W$ double-well potential having wells located at $0$ and $1$ with $W(0)=W(1)$.
Then, equation \eqref{scalar_damp_rd_1d} supports monotone steady states connecting
equilibria $\phi_-=1$ and $\phi_+=0$.
\end{prop}

\begin{proof}
We report here the standard proof for reader's convenience.
Substituting $c=0$, equation \eqref{twode} reduces to
\begin{equation*}
	a\frac{d^2\phi}{d\xi^2} - W'(\phi) = 0,
\end{equation*}
Multiplying by the derivative $d\phi/d\xi$, we end up with the conservative form
\begin{equation*}\label{twode_c=0}
	\frac{d}{d\xi}\left\{\frac{a}{2}\left(\frac{d\phi}{d\xi}\right)^2-W(\phi)\right\} = 0,
\end{equation*}
which can be integrated.
Then, we infer
\begin{equation}\label{twode_c=0_ter}
		\frac{d\phi}{d\xi} = - \sqrt{{2}/{a}}\,\cdot\,\sqrt{W(\phi) - W(\phi_\pm)},
\end{equation}
recalling that $\phi$ is monotone decreasing since $\phi_+ =0< 1=\phi_-$.
Hence, among other solutions, equation \eqref{twode_c=0_ter} defines implicitly the steady profile $\phi=\phi(\xi)$ by
\begin{equation*}
	\int_{\phi(\xi_0)}^{\phi(\xi)}\frac{ds}{\sqrt{W(s) - W(\phi_\pm)}}=\sqrt{{2}/{a}}\,(\xi_0-\xi)
\end{equation*}
connecting $\phi_- =1$ to $\phi_+ =0$ for any given $\xi_0\in\mathbb{R}$.
\qed
\end{proof}

As an example, let us consider the case $f(u)=\kappa\,u(u-1/2)(1-u)$.
Since the potential is given by $W(u)=\tfrac14\,\kappa\,u^2(1-u)^2$, there holds
\begin{equation*}
	\int_{1/2}^{\phi(x)}\frac{2\,ds}{s(1-s)}=\sqrt{{2\kappa}/{a}}\,(x_0-x)
\end{equation*}
that gives
\begin{equation*}
	\ln\left(\frac{\phi(x)}{1-\phi(x)}\right)=\sqrt{{\kappa}/{2a}}\,(x_0-x).
\end{equation*}
Expliciting the value $\phi=\phi(x)$, we obtain
\begin{equation}\label{explicitprofile}
	\phi(x)=\frac{1}{1+e^{\sqrt{\frac{\kappa}{2\,a}}(x-x_0)}}
	\qquad (x_0\in\mathbb{R}).
\end{equation}
As stated at the beginning, the propagation speed is $c=0$.

\subparagraph{Third-order polynomial reaction function}

Next, we focus on the case $c_\ast>0$, which occurs, again by formula \eqref{speedformula1}  for $W(1)<W(0)$.
In the case  of the third order polynomial
\begin{equation}\label{polyreact}
	f(u)=\kappa\,u(u-\alpha)(1-u)
\end{equation}
with $\kappa>0$, this translates into the choice $\alpha\in(0,1/2)$.

To start with, let us focus on the limiting case $\sigma=\tau=0$, that is on the parabolic reaction-diffusion equation
\begin{equation}\label{allencahn}
	 \partial_t u = a\,\partial_{xx} u + \kappa u(u-\alpha)(1-u).
\end{equation}
In such a case, there exist explicit formulas for both propagation speed $c$ and front profile $\phi$.
Indeed, let us set
\begin{equation*}
	\frac{d\phi}{d\xi}=-A\phi(1-\phi).
\end{equation*}
for some constant $A>0$.
Since
\begin{equation*}
	\frac{d^2\phi}{d\xi^2}=-A(1-2\phi)\frac{d\phi}{d\xi}=A^2\phi(1-\phi)(1-2\phi),
\end{equation*}
inserting in \eqref{twode} with $\sigma=\tau=0$ and simplifying the factor $\phi(1-\phi)$, we infer
\begin{equation*}
	(\kappa-2aA^2)\phi+aA^2-cA-\kappa\alpha=0
\end{equation*}
which gives $A=\sqrt{\kappa/(2a)}$ and 
\begin{equation}\label{explicitspeed}
	c=c_0:=\sqrt{2\,a\kappa}\left(\tfrac{1}{2}-\alpha\right).
\end{equation}
Thus, the corresponding profile $\phi$ solves the Bernoulli equation ${d\phi}/{d\eta}=-\phi+\phi^2$
where $\eta=(\kappa/2a)^{1/2}\,\xi$, which is explictly given by
\begin{equation*}
	\phi(\xi)=\frac{1}{1+e^{\sqrt{\frac{\kappa}{2\,a}}(\xi-\xi_0)}}\qquad (\xi_0\in\mathbb{R}),
\end{equation*}
which, incidentally, coincide with \eqref{explicitprofile} when $\xi=x$.

When dealing with propagation fronts for \eqref{scalar_damp_rd_1d} with $\sigma=0$, that is
\begin{equation*}\label{purelydamped}
	\tau\,\partial_{tt} u + \partial_t u = a\,\partial_{xx} u+f(u),
\end{equation*}
a formula, corresponding to \eqref{explicitspeed}, can be provided.
Indeed, equation \eqref{twode} with $\sigma=0$ coincide with the traveling wave equation
for \eqref{allencahn} where $a$ has been replaced by $a-\tau c^2$.
Thus, adding the subscript $\tau$ to $c$ to give evidence to dependency, there holds
\begin{equation*}
	c_\tau=\sqrt{2\,(a-\tau c_\tau^2)\kappa}\,\cdot\,\left(\tfrac12-\alpha\right).
\end{equation*}
Squaring and rearranging, we infer
\begin{equation*}
	\left\{1+2\kappa\tau\left(\tfrac12-\alpha\right)^2\right\}c_\tau^2=2\,a\kappa\left(\tfrac12-\alpha\right)^2,
\end{equation*}
and thus
\begin{equation}\label{explicitspeed2}
	c_\tau=\frac{c_0}{{\sqrt{1+\tau c_0^2/a^2}}}.
\end{equation}
where $c_0$ is given in \eqref{explicitspeed}.
There is a strict connection between relation \eqref{explicitspeed2} and \eqref{rescaledspeed},
being one the inverse of the other in the case $a=1$.
Specifically, relation \eqref{explicitspeed2} goes beyond the special case of the cubic $f$,
holding for general reaction function.
In particular, since $0\leq c_\tau<c_0$ for $\tau>0$, as shown by the inequality
\begin{equation*}
	\frac{c_\tau-c_0}{c_0}=\frac{1}{{\sqrt{1+\tau c_0^2/a^2}}}-1<0,
\end{equation*}
the propagation phenomena is always slowed down when pure damping is added,
inertia being limited to the deceleration effect of the front.

When dealing with hyperbolic reaction-diffusion equation \eqref{scalar_damp_rd_1d} with $\sigma\in(0,\tau]$ and cubic $f$,
to our knowledge, there is no available extension of the explicit formulas \eqref{explicitspeed} and \eqref{explicitspeed2}.
In particular, as it will be shown later on, the addition of the relaxation term, i.e. $\sigma=\tau$, the situation
 relative to the difference in propagation speed can change in some regime of the parameter $\alpha\in(0,1)$.

\subparagraph{Piecewise affine reaction function with a bistable shape}

\begin{figure}[b]
\sidecaption
\includegraphics[width=7.0cm]{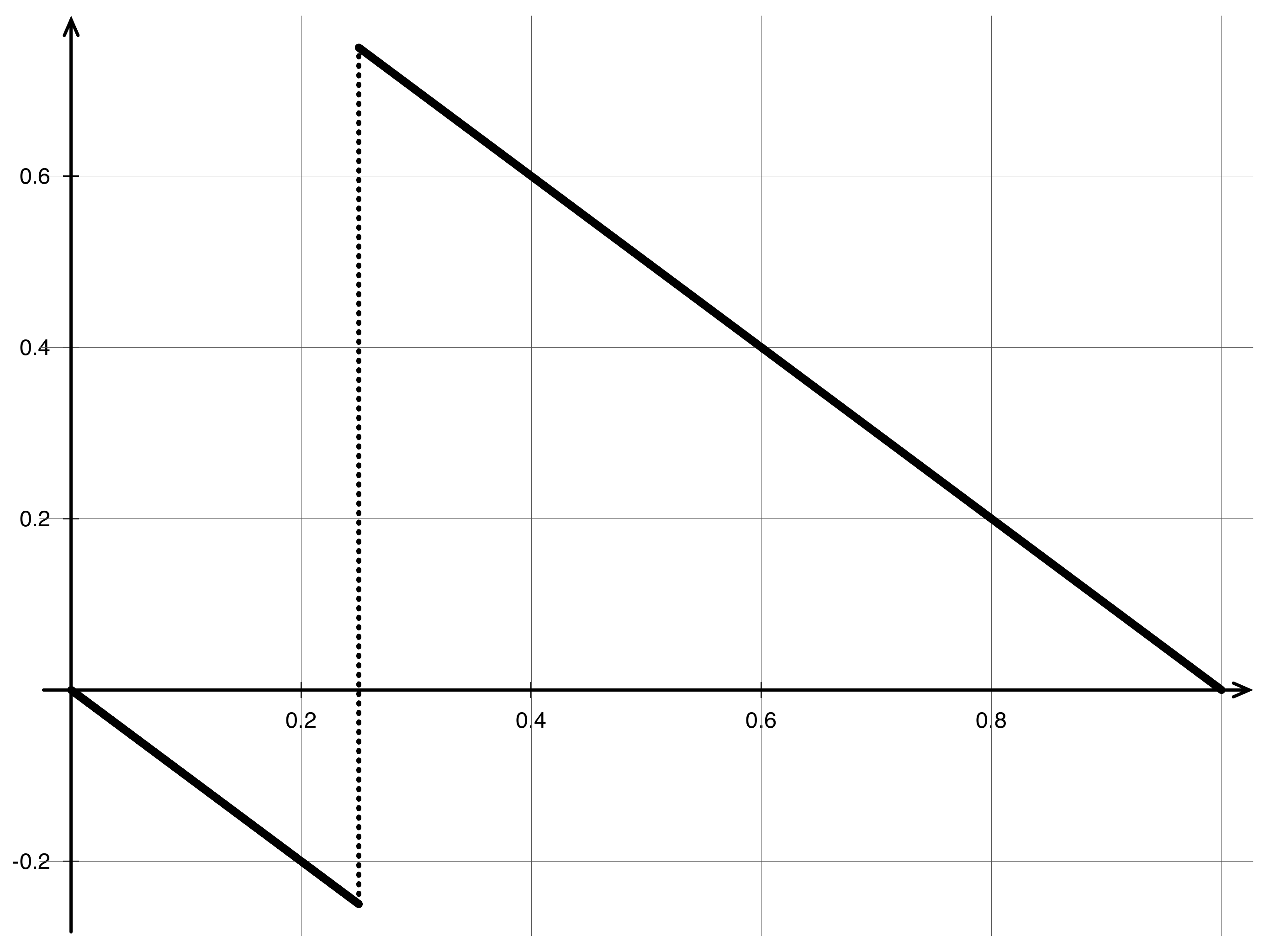}
\caption{Graph of the function $f$ given in \eqref{pwl} with parameters $m=1$ and $\alpha=0.25$.}
\label{fig:pwl}
\end{figure}

Finally, following the approach in \cite{McKe70}, we compute explicit traveling wave solutions
for a very specific form for the reaction function $f$ of bistable type.
Specifically, we concentrate on a piecewise affine function given by
\begin{equation}\label{pwl}
	f(u)=\left\{\begin{aligned}
		&-m\,u	&\quad	&u<\alpha,\\
		&m(1-u)	&\quad	&u\geq \alpha.
		\end{aligned}\right.
		\qquad m>0,\, \alpha\in(0,1),
\end{equation}
(see Fig.\ref{fig:pwl}).
In such a special case, it is possible to provide an explicit expression for both the traveling wave profile $(\phi,\psi)$
and of its speed $c$  also for the hyperbolic model \eqref{twode}.
 Indeed, let us go back to \eqref{twode_1st} and rewrite it as
\begin{equation*}\label{twode_pwl_left}
	\frac{d\phi}{d\xi}=\psi,\qquad 
	(a-\tau c^2)\,\frac{d\psi}{d\xi}= m\,\phi - c\left(1+\sigma m\right)\psi ,
\end{equation*}
to be matched at $\phi=\alpha$ with
\begin{equation*}\label{twode_pwl_right}
	\frac{d\phi}{d\xi}=\psi,\qquad 
	(a-\tau c^2)\,\frac{d\psi}{d\xi}= m(\phi-1) - c\left(1+\sigma m\right)\psi.
\end{equation*}
Since the two singular points are saddles, the matching amounts in choosing the critical value of the parameter $c$
such that the unstable manifold of the singular point $(0,0)$ intersects, at $\phi=\alpha$, the stable manifold of $(1,0)$.

The directions of the unstable/stable manifolds are described by the eigenvectors of the corresponding linearized equation.
Hence, denoted by $(\tilde \phi,\tilde \psi)$ the perturbation of the equilbrium state $(\bar\phi,0)$, they are given
by the eigendirection of the matrix
\begin{equation*}
	\mathbf{A}:=\frac{1}{a-\tau c^2}\begin{pmatrix} 0 & a-\tau c^2 \\ m & -c(1+\sigma m) \end{pmatrix}
\end{equation*}
In particular, this means that $(\phi,\psi)$ belongs to the unstable/stable manifold if and only if $\tilde \psi=\lambda_\pm \tilde \phi$,
where $\lambda_\pm$ denote the (positive/negative) roots of the characteristic polynomial
\begin{equation*}
	p(\lambda)=\det(\mathbf{A}-\lambda\,\mathbf{I})=\frac{1}{a-\tau c^2}\left\{(a-\tau c^2)\lambda^2+c(1+\sigma m)\lambda-m\right\}.
\end{equation*}
Specifically, the explicit values for $\lambda_\pm$ are
\begin{equation*}
	p(\lambda_\pm)=0\qquad\iff\qquad
	\lambda=\lambda_\pm:=\frac{-c(1+\sigma m)\pm\sqrt{\Delta(c)}}{2(a-\tau c^2)}
\end{equation*}
where the discriminant $\Delta$ is
\begin{equation*}
	\begin{aligned}
	\Delta(c)&:=c^2(1+\sigma m)^2+4(a-\tau c^2)m\\
		&=\left[(1-\sigma m)^2-4(\tau-\sigma)m\right]c^2+4am,
	\end{aligned}
\end{equation*}
which is strictly positive in the regime $c^2<a/\tau$.
Thus, the stable manifold of $(0,0)$ and the unstable manifold at $(1,0)$ are given by 
$\tilde \psi=\lambda_+\tilde \phi$ and $\tilde \psi=\lambda_-\tilde \phi$, that is
\begin{equation*}
	\psi=\lambda_- \phi\qquad\textrm{and}\qquad \psi=\lambda_+(\phi-1).
\end{equation*}
The two graphs intersect at $\phi=\alpha$ if and only if $|\lambda_-|\alpha =\lambda_+(1-\alpha)$.
Recalling the explicit formulas for $\lambda_-$ and $\lambda_+$, the latter equality can be rewritten as
\begin{equation*}
	\sqrt{\Delta(c_{\textrm{ex}})}(1-2\alpha)=c_{\textrm{ex}}(1+m\sigma).
\end{equation*}
After some straightforward algebraic manipulations, we end up with
\begin{equation}\label{exactPWL}
	c_{\textrm{ex}}=\left\{\frac{ma}{(1+m\sigma)^2\alpha(1-\alpha)+m\tau(2\alpha-1)^2}\right\}^{1/2}(1-2\alpha).
\end{equation}
\begin{figure}[htb]\centering
\includegraphics[width=9cm]{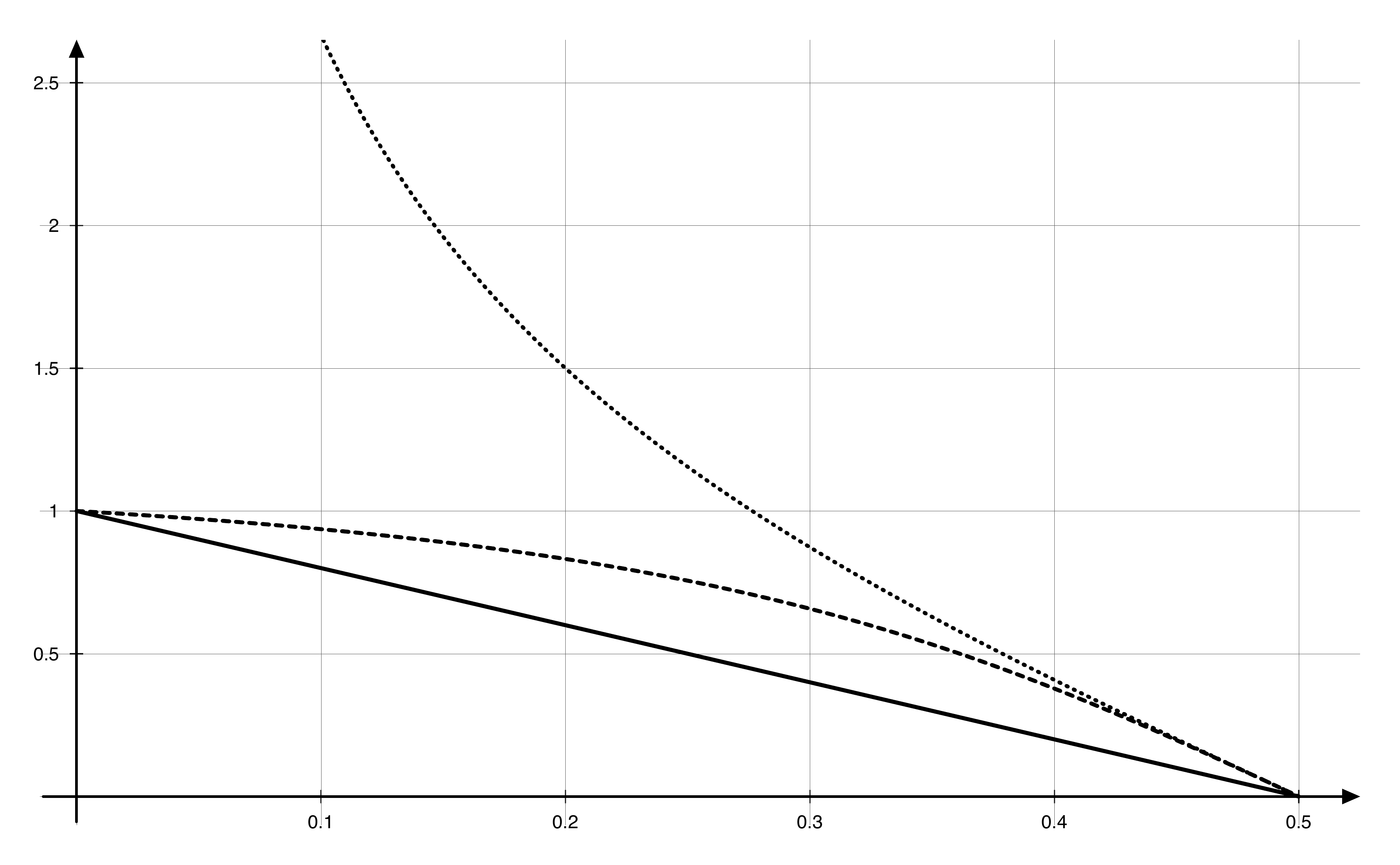}
\caption{Exact value of the speed, as in \eqref{exactPWL}, with $a=m=1$ and $\sigma=\tau=0$ (dotted),
$\sigma=0$, $\tau=1$ (dashed), $\sigma=\tau=1$ (continuous).}
\label{fig:alphaSpeedPWL}
\end{figure}

Comparing the speeds $c_{\textrm{ex}}$ for a generic choice of  parameters $\sigma$ and $\tau$ and $c_0$ for $\sigma=\tau=0$ gives
\begin{equation*}
	\frac{c_{\textrm{ex}}}{c_0}=\left\{\frac{\alpha(1-\alpha)}{(1+\sigma m)^2\alpha(1-\alpha)+\tau m(2\alpha-1)^2}\right\}^{1/2}
\end{equation*}
For $\sigma\in[0,\tau]$, since $\alpha(1-\alpha)<1/4$ for $\alpha\neq 1/2$, there holds
\begin{equation*}
	\frac{\alpha(1-\alpha)}{(1+m\sigma)^2\alpha(1-\alpha)+m\tau(2\alpha-1)^2}
		<\frac{1}{(1+m\sigma)^2+4m\tau(2\alpha-1)^2}\leq 1
\end{equation*}
with the equality holding if and only if $\tau=0$.
Hence, in the same regime, it follows
\begin{equation*}
	\frac{c_{\textrm{ex}}-c_0}{c_0}=\left\{\frac{\alpha(1-\alpha)}{(1+m\sigma)^2\alpha(1-\alpha)+m\tau(2\alpha-1)^2}\right\}^{1/2}-1<0.
\end{equation*}
In particular, the (hyperbolic) propagation speed $c_{\textrm{ex}}$ is always smaller than the corresponding (parabolic) speed $c_0$
for any choice of the couple $\sigma$ and $\tau$.
This could be also recognised, observing directly that the value of $c_{\textrm{ex}}$, regarded as a function of $\sigma$ and $\tau$,
is strictly decreasing with respect to both variables.

Let us remark that, in such a case, the function $f$ is discontinuous (increasing) at the value $u=\alpha$
and, thus, the first derivative of $f$ is, lousely speaking, equal to $+\infty$.
In particular, the dissipativity condition $1-\sigma f'>0$ is never satisfied at such a point whenever $\sigma>0$,
with dramatic consequences to be explored in the next Section.

\section{Numerical computation of the propagation speed}
\label{sect:numerical}

\setcounter{equation}{0}
\renewcommand{\theequation}{3.\arabic{equation}}

From now on, we restrict the attention to two main cases corresponding to the choices: $\sigma=0$, $\tau>0$ and $\sigma=\tau>0$,
reported here for reader's convenience,
\begin{equation*}
	\begin{aligned}
	\tau\partial_{tt} u + \partial_t u & = a \partial_{xx} u + f(u)
		&\qquad &\textrm{(damping)}\\
	\tau\partial_{tt} u + \partial_t \left\{u-\tau f(u)\right\} & = a \partial_{xx} u + f(u)
		&\qquad &\textrm{(relaxation)}
	\end{aligned}
\end{equation*}
where $f(u)=\kappa\,u(u-\alpha)(1-u)$ with $\kappa>0$ and $\alpha\in(0,1)$.
Coherently with the previous part of the paper, we focus on propagating waves connecting 1 at $-\infty$ with 0
at $+\infty$ in the case $\alpha\in(0,1/2]$, so that the speed $c_{\textrm{ex}}$ is non-negative as a consequence
of the relation $W(1)\leq W(0)$, see identity \eqref{speedformula1}.

\subsection{Computation of the propagation speed}
\label{subsect:comp_prop_speed}

In the purely damped case, the explicit formula  \eqref{explicitspeed2} for the propagation speed can be used to assess
the reliability of the so-called {\it phase-plane algorithm}, presented in detail in the next subsection. 
On the other hand, when relaxation is taken into account, there is no explicit formula for the velocity.
Thus, an approximated version of its value should be considered as furnished by some algorithm.
Based on the tests used in the damped case, we will consider as ``exact'' speed $c_{\textrm{ex}}$
the ones provided by the phase-plane algorithm  (later on, denoted by $c_{\textrm{du},\theta}$),
and use it to test the capability of two (dynamical) numerical schemes to provide  genuine predictions.

\subparagraph{Phase plane algorithm}

 As stated before, both singular points of the ODE system for traveling waves \eqref{twode_system} 
are saddles  in the bistable case.
As a consequence, both the corresponding unstable/stable manifold are one-dimensional.
Therefore, the existence of a heteroclinic connection is equivalent to the fact that, for an appropriately tuned parameter
$c=c_{\textrm{ex}}$, the unstable curve exiting from the critical point $(1,0)$ intersects the stable curve
entering the critical point $(0,0)$.
Based on the rotated vector field property, we can perform a shooting-type argument 
and transform the problem of the existence of a heteroclinic  orbit into the search of a zero of a given function.
Such a step  can be performed by preliminarily finding a reliable approximation of the solution to
an ordinary differential equation and then by means of a standard interval division scheme, furnishing the exact value $c_{\textrm{ex}}$
of the propagation speed.

To enter the details, we denote by $v_0=v_0(\phi,c)$, the stable manifold of $(0,0)$ and by $v_1=v_1(\phi,c)$ the unstable manifold of $(1,0)$.
Then, we look for two different solutions of the first order equation
\begin{equation}\label{firstorderTW}
	\frac{\partial v}{\partial \phi}=\frac{d\psi/d\xi}{d\phi/d\xi}
		=\frac{1}{a-\tau c^2}\left\{\frac{1}{\psi}\frac{dW}{du}(\phi)-c\left[1+\sigma\,\frac{d^2 W}{du^2}(\phi)\right]\right\}
\end{equation}
with initial conditions along the stable/unstable manifold of $(0,0)$/$(0,1)$.

Curves $v_0$ and $v_1$ are determined by choosing an initial datum
on the corresponding stable/unstable manifold as provided by the linearized operator at the two critical points.
Namely, at $\bar u$, we compute the eigenvectors relative to the eigenvalues $\lambda_\pm=\lambda_\pm(\bar u;c)$
as given by \eqref{proots}.
Then, we  approximate the solutions $v_0=v_0(\cdot,c)$ and $v_1=v_1(\cdot,c)$
 with the ones defined by the initial data
\begin{equation*}
	v_0(\varepsilon,c)=\lambda_-(0,c)\varepsilon
	\qquad\textrm{and}\qquad
	v_1(1-\varepsilon,c)=-\lambda_+(1,c)(1-\varepsilon)
\end{equation*}
for $\theta$ small enough and solving forward/backward \eqref{firstorderTW} for $v_0$/$v_1$, respectively.

 Denoting by $v_0$ and $v_1$ such approximations,
we evaluate the difference function $h$ of $v_0$ and $v_1$ at $u=\alpha$,  that is
\begin{equation*}
	h(c):=v_0(\alpha,c)-v_1(\alpha,c),
\end{equation*}
for $c\in(-\sqrt{{a}/{\tau}},\sqrt{{a}/{\tau}})$.
It can be readily seen that
\begin{equation*}
	h(-\sqrt{a/\tau})<0<h(\sqrt{a/\tau}).
\end{equation*}
Moreover,  relying on the rotated vector field property, the function $h$ is strictly increasing  in $(-\sqrt{{a}/{\tau}},\sqrt{{a}/{\tau}})$
and, thus, it has a single zero, corresponding to the value $c_{\textrm{ex}}$.
The heteroclinic orbit corresponds to such a choice of the critical speed $c_{\textrm{ex}}$ such that $h(c_{\textrm{ex}})=0$,
which is uniquely determined since the function $h$ is strictly monotone increasing,

\subparagraph{Heuristic validation of the phase-plane algorithm in the purely damped case}

\begin{figure}[htb]\centering
\includegraphics[width=9cm]{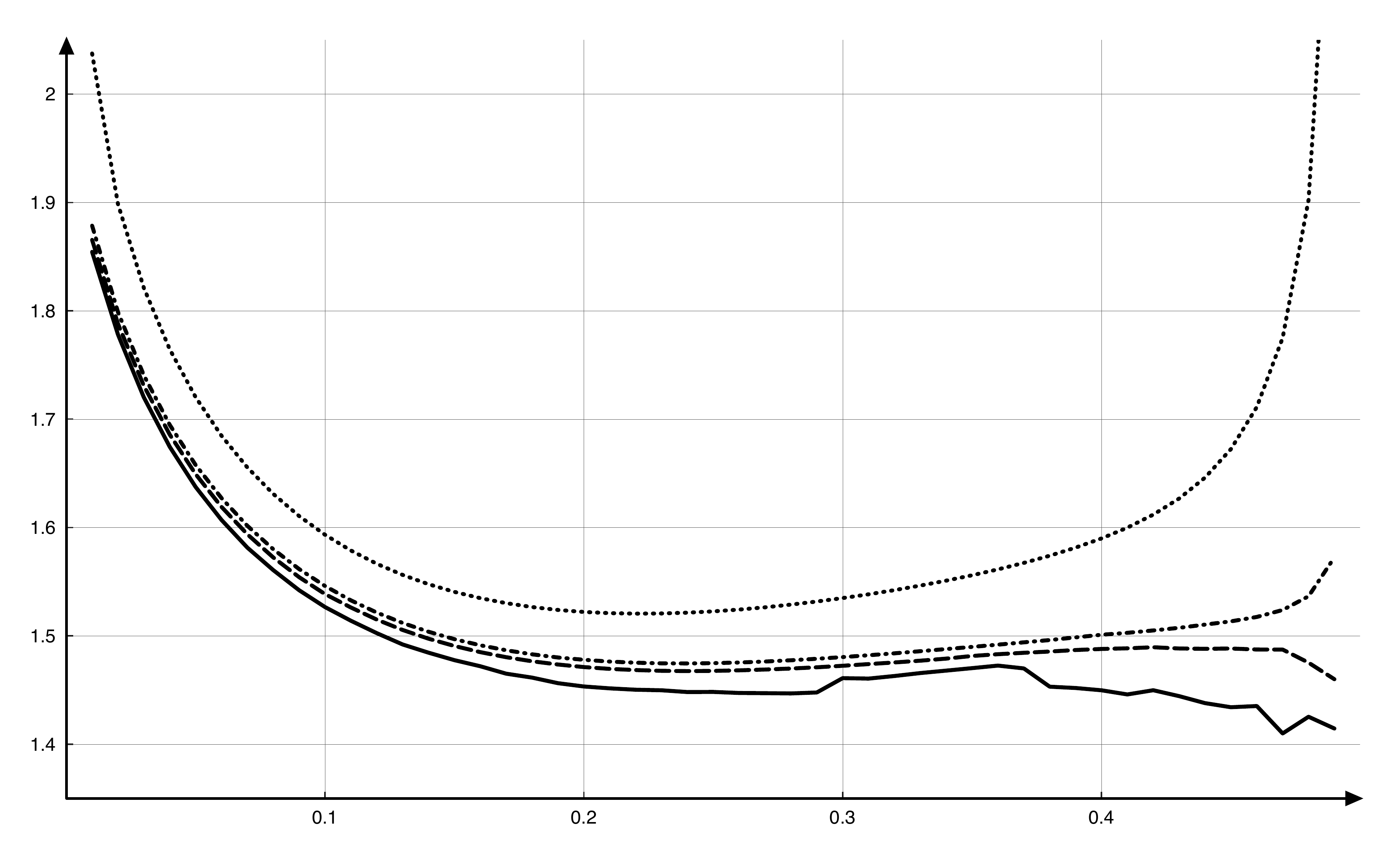}
\caption{Case $\tau=a=\kappa=1$: graphs of the values of $E_{\textrm{du},\varepsilon}/\textrm{du}$ as a function of $\alpha\in(0,0.5)$ where
the relative error $E_{\textrm{du},\theta}$ is given in \eqref{relative_error_damped} for $\varepsilon=10^{-8}$
and discretization step equal to different choices of $\textrm{du}$:
$10^{-2}$ (dotted), $10^{-3}$ (dotted-dashed), $10^{-4}$ (dashed), $10^{-5}$ (continuous).}
\label{fig:plotErrRelDamp}
\end{figure}

Next, we compare the exact formula \eqref{explicitspeed2} in the case $\sigma=0$, $a=\kappa=1$,
recalled here for reader's convenience, viz.
\begin{equation*}\label{exactspeed_damped}
	c_{\textrm{ex}}=\frac{\sqrt{2}\,\left(1/2-\alpha\right)}{{\sqrt{1+2\tau\left(1/2-\alpha\right)^2}}},
\end{equation*}
with the approximated value $c_{\textrm{du},\varepsilon}$ provided by the 
phase-plane algorithm using the value $E_{\textrm{du},\varepsilon}$ as measure of the relative error, defined by
\begin{equation}\label{relative_error_damped}
	E_{\textrm{du},\varepsilon}:=\left|\frac{c_{\textrm{du},\varepsilon}-c_{\textrm{ex}}}{c_{\textrm{ex}}}\right|
\end{equation}
To start with, we learn from Fig.\ref{fig:plotErrRelDamp} that there is numerical evidence of a scheme of order $1$
in the case $\tau=1$.
Different values of $\tau$, $a$ and $\kappa$ fits into the same scenery.

From this, we extrapolate the final (reliable) choices $\textrm{du}=10^{-5}$ and $\varepsilon=10^{-8}$.
The corresponding values for the exact formula $c_{\textrm{ex}}$, the approximated value $c_{\textrm{du},\varepsilon}$
and the relative error $E_{\textrm{du},\varepsilon}$, are reported in Table \ref{tab:du10-5}, for different values of 
the unstable zero $\alpha$, chosen as a value in $(0,1/2)$.

\begin{table}
\caption{Case $\tau=a=\kappa=1$: values for  $c_{\textrm{ex}}$, $c_{\textrm{du},\varepsilon}$ and $E_{\textrm{du},\varepsilon}$
relative to nine different choices of the unstable zero $\alpha$ relative to the choices
$\textrm{du}=10^{-5}$ and $\theta=10^{-8}$.}
\begin{center}\begin{tabular}{p{1.50cm}p{2.00cm}p{2.00cm}p{1.50cm}}
$\alpha$	& $c_{\textrm{ex}}$	& $c_{\textrm{du},\varepsilon}$	& $E_{\textrm{du},\varepsilon}$		\\
\noalign{\smallskip}\svhline\noalign{\smallskip}
0.05		& 0.5368950		& 0.5369038				& $1.64\times 10^{-5}$			\\
0.10		& 0.4923660		& 0.4436135				& $1.53\times 10^{-5}$			\\
0.15		& 0.4436070		& 0.4436135				& $1.48\times 10^{-5}$			\\
0.20		& 0.3905667		& 0.3905724				& $1.45\times 10^{-5}$			\\
0.25		& 0.3333333		& 0.3333382				& $1.45\times 10^{-5}$			\\
0.30		& 0.2721655		& 0.2721695				& $1.46\times 10^{-5}$			\\
0.35		& 0.2075143		& 0.2075174				& $1.47\times 10^{-5}$			\\
0.40		& 0.1400280		& 0.1400300				& $1.45\times 10^{-5}$			\\
0.45		& 0.0705346		& 0.0705356				& $1.43\times 10^{-5}$			\\
\noalign{\smallskip}\hline\noalign{\smallskip}
\end{tabular}\end{center}
\label{tab:du10-5}
\end{table}

In the case $\sigma\in(0,\tau]$ for some $\tau>0$, to our knowledge, there is no explicit formula for the case
of the double-well potential $W$,  given by \eqref{verhulst_pot}.
Hence, we  consider the speed approximation provided by the phase-plane algorithm
with the values for $\textrm{du}$ and $\varepsilon$ previously detected.
From now on, for simplicity, we will denote $c_{\textrm{du},\varepsilon}$ by $c_{\textrm{ex}}$ and consider the relative errors
with respect to such an approximated value.

\begin{figure}[htb]\centering
\includegraphics[width=9cm]{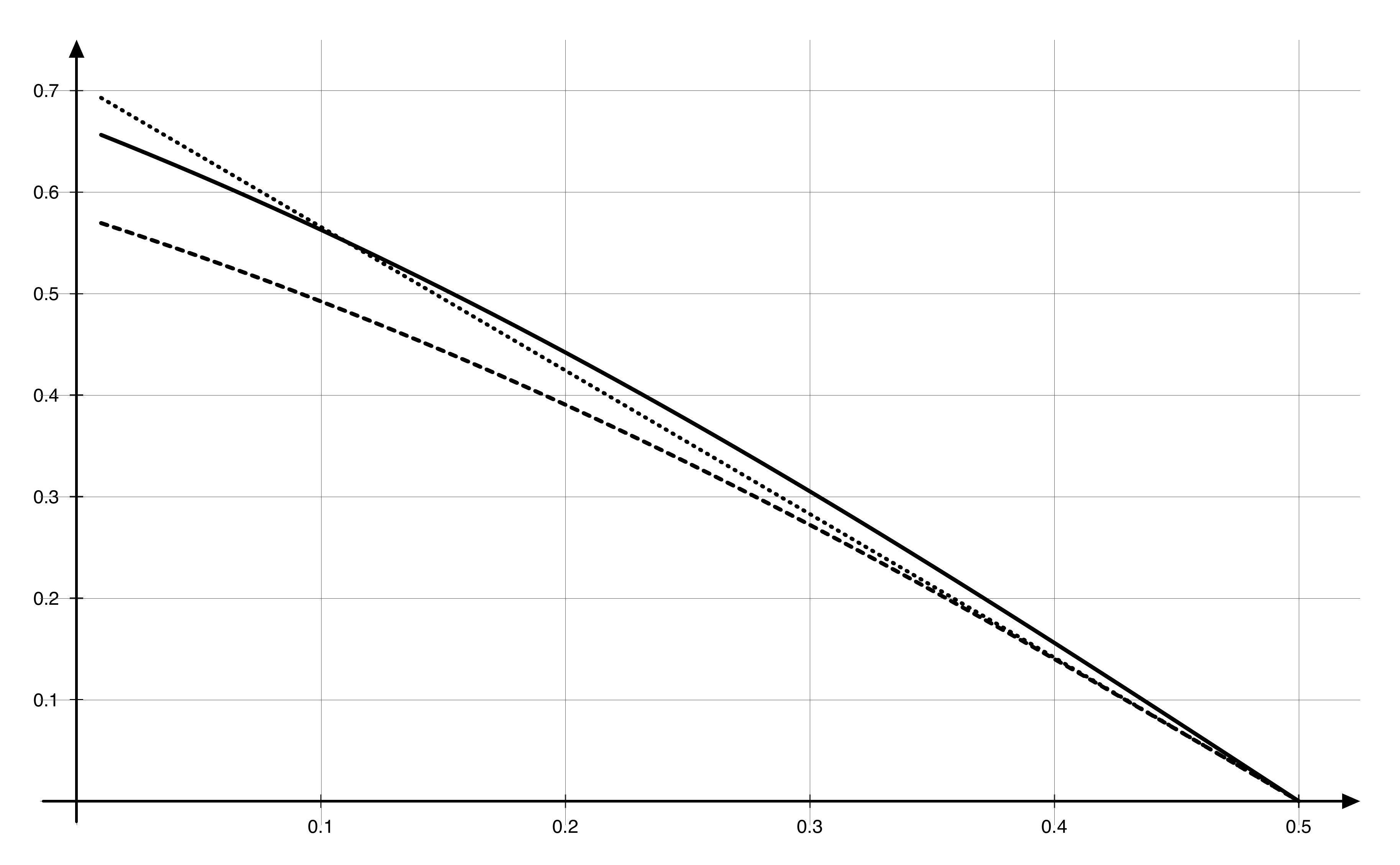}
\caption{Case $\tau=a=\kappa=1$: comparison of the graphs of the speeds:
parabolic Allen--Cahn (dotted), see \eqref{explicitspeed};
hyperbolic Allen--Cahn with damping (dashed), see \eqref{explicitspeed2};
hyperbolic Allen--Cahn with relaxation (continuous).}
\label{fig:alphaSpeed}
\end{figure}

To conclude, in Figure \ref{fig:alphaSpeed}, we compare the values for the Allen--Cahn equation in the standard parabolic case, 
in the hyperbolic case with damping, in the hyperbolic case with relaxation.
It is transparent that the role played in the latter is crucially different and it exhibits values $\alpha$ where the role
of inertia is purely dissipative and others values for which sustained propagation is present.

\subsection{PDE-based algorithms to approximate the propagation speed}
\label{subsect:PDEalgo}

The aim of this Subsection is to compare the capability of two different PDE-based algorithms to recover 
 a reliable approximation of the speed of a front.
The strategy is different with respect to the one presented in Subsection \ref{subsect:comp_prop_speed} being of {\it dynamical} nature,
i.e. grounded on the preliminary  determination of the numerical solution of the underlying  partial differential equation.
Entering the details, we choose a scheme for the PDE and solve it in the space interval $[0,L]$, with zero-flux boundary conditions,
in the time span $[0,T]$, corresponding to some initial datum.
Then, choosing two consecutive frames $u(\cdot,s)$ and $u(\cdot,t)$ with $0<s<t$, we look for a strategy furnishing a scalar
value $c$ such that 
\begin{equation*}
	u(x,t)-u(y,s)\approx \phi(x-ct)-\phi(y-cs).
\end{equation*}
The key point stems in reducing from two functions (i.e. the solution profiles) to a single scalar value which should be able to describe,
in principle, the overall propagating characteristic of the  wave.

We consider the three  numerical schemes described in Subsection \ref{subsec:algorithms}
(with the kinetic algorithm limited to the relaxation case), freezing the data relative to the two profiles $u(\cdot,s)$
and $u(\cdot,t)$ with $0<s<t$ appropriately chosen.
Then, we determine an approximation of the speed by means of some  suitably chosen algorithm.

Two main tools can be used to provide an  estimate of the speed,
the {\it scout \& spot algorithm} and the {\it LeVeque--Yee formula}, which we present in details in the following paragraphs.

At this point, a word of caution is required.
Indeed, the approximated expression $c=c_{\mathrm{dx},\mathrm{dt}}$ for the velocity is relative to the specific
numerical scheme and, in addition to the scheme itself, it depends on both choices of space and time mesh sizes.
Also, the potential existence of a propagating front for the semi- and fully-discrete schemes (not explored in this Chapter)
is not necessarily related to the existence of a continuous propagating front (sketched in this Chapter and rigorously proved
in \cite{BouiCalvNadi14, Hade88, LattMascPlazSime19,LattMascPlazXX} for different types of hyperbolic
reaction-diffusion equations).
Results on the existence of parabolic reaction-diffusion traveling waves can be found in \cite{BateChenChma03, Keen87, Mall99, Zinn92}
for spatially-discrete schemes (sometimes referred to as ``lattices'') and in \cite{ChowMallShen98, ElmeVanV05, HupkEtAl20, HupkVanV16}
for the fully-discrete case.
Moreover, discussions relative to hyperbolic equations can be found in \cite{CarpDuro05,ElmeVanV99}.
For completeness, let us also mention that a corresponding exploration in the context of conservation laws,
started in \cite{Jenn74, MajdRals79}, can be found in \cite{Benz98, Serr07}.
To our knowledge, a detailed scrutiny of existence of propagating front for both semi- and fully-discrete 
schemes proposed in Subsection \ref{subsec:algorithms} is not currently available and we regard at it
as a very interesting issue.
 
In any case, as shown in most of the previous references, it is reasonable to associate to any convergent numerical
scheme a number --coinciding with the ``discrete'' speed of propagation-- that could be regarded as an approximation
of the exact velocity value in the continuous setting under appropriate limiting behavior of the parabolic
ratio $\textrm{dt}/\textrm{dx}^2$.
 
Precisely, given one of the three numerical schemes to approximate the hyperbolic reaction-diffusion equation (first-order, Li\'enard, kinetic)
together with one of the two possible algorithms to estimate the speed (scout\&spot, LeVeque--Yee, with details provided in the subsequent pages),
we consider as a reliable error measure the value
\begin{equation*}
	E_{\ast}^{\textrm{scheme}}:=\left|\frac{c_{\ast}^{\textrm{scheme}}-c_{\textrm{ex}}}{c_{\textrm{ex}}}\right|,
\end{equation*}
where, as stated before, $c_{\textrm{ex}}$ coincides with $c_{\textrm{du},\theta}$ with $\textrm{du}=10^{-5}$ and $\varepsilon=10^{-8}$
and $c_{\ast}^{\textrm{scheme}}$ is the estimated value for the propagation speed.
We anticipate that we are going to compare the three schemes considering spatial and temporal mesh size given, respectively,
by $\mathrm{dx}=10^{-1}$ and $\mathrm{dt}=10^{-3}$, so that the ratio $\textrm{dt}/\textrm{dx}^2$ has the exact value $10^{-1}$
to be regarded as  a ``small number''.

\subparagraph{Scout \& spot algorithm}
The first determines the speed of propagation considering a fixed level curve, say $\theta$, taking into account the fact that, whenever
the solution $u$ converges to the propagating front $\phi$, the relation $u(x,t)\approx \phi(x-ct)$ holds asymptotically in time, i.e. as $t\to+\infty$.
Let $\phi_-<\phi_+$ and fix a value $\theta\in(\phi_-,\phi_+)$ and consider two different
time instants, denoted here by $t$ and $s$, such that $u(x(s),s)=u(x(t),t)=\theta$, then
\begin{equation*}
	x(t)-ct\approx \phi^{-1}(\theta)\approx x(s)-cs.
\end{equation*}
Hence, we deduce the approximation formula
\begin{equation}\label{approxcurve}
	c\approx \frac{x(t)-x(s)}{t-s}.
\end{equation}
Translating such approximated rule in a definite algorithm is based on the introduction of a specific space mesh 
$J=\{x_1,\dots,x_j\}=\{\textrm{dx},2\textrm{dx},\dots,j\textrm{dx}\}$.
Assuming that the profile $u^n_j$ is strictly monotone increasing with respect to $j$, the first step consists in considering
the first value where the threshold $\theta$ is trespassed for any given time $t^n$, that is
\begin{equation*}
	n\;\longmapsto\; j^n(\theta):=\max\{j\in J\,:\,u^n_j<\theta\}.
\end{equation*}
Approximation formula \eqref{approxcurve} becomes 
\begin{equation}\label{levelcurvealgo}
	\begin{aligned}
	c_{\textrm{s\&s}}^{n,p} = c_{\textrm{s\&s}}^{n,p} (\theta)
		&= \frac{j^{n+p}(\theta)-j^n(\theta)}{t^{n+p}-t^n}\,\cdot\,\textrm{dx}\\
		&= \left[j^{n+p}(\theta)-j^n(\theta)\right]\,\cdot\,\frac{\textrm{dx}}{p\,\textrm{dt}}
	\end{aligned}
\end{equation}
Such procedure corresponds to a {\it piecewise constant interpolation} of the states $u^n_j$ and $u^n_{j+1}$.
Moreover, the above formula shows that the propagation speed of slow waves provided by such a level curve algorithm
is ``quantized'', that is any candidate as limiting speed is an {\it integer} multiple of the positive value by $\textrm{dx}/(p\,\textrm{dt})$.

Applying such an algorithm requires a number of choices, which can be matter of criticism, starting from the fact that the profile is
expected to be monotone increasing.
Here, we choose $\theta=\alpha$, $p=T/(2\textrm{dt})$ so that the speed is approximated up to an error of order
$\textrm{dx}/(p\,\textrm{dt})=10^{-2}$ in the case $T=50$ and $\textrm{dx}=10^{-1}$.

\subparagraph{LeVeque--Yee formula}

The second strategy, inspired by \cite{LeVYee90}, makes use of a spatial average of the profile and it does not require
any monotone assumption on the solution.
Anyway, it is still needed that the two asymptotic states, $\phi_-$ at $-\infty$ and $\phi_+$ at $+\infty$, are different,
i.e. the connection has to be heteroclinic.

Let $\phi$ be a differentiable function with asymptotic states $\phi(\pm\infty)=\phi_\pm$.
The LeVeque--Yee formula takes advantage from the exact relation 
\begin{equation*}\label{increment}
	\int_{\mathbb{R}} \left\{\phi(x+h)-\phi(x)\right\}\,dx=h\,\left[\phi\right]
\end{equation*}
where $[\phi]:=\phi_+-\phi_-$.
The above formula can be proved by observing that
\begin{equation*}
	\begin{aligned}
	\int_{\mathbb{R}} \left\{\phi(x+h)-\phi(x)\right\}\,dx&=h\int_{\mathbb{R}} \int_0^1 \frac{d\phi}{dx}(x+\theta h)\,d\theta\,dx\\
		&=h\int_0^1 \int_{\mathbb{R}}  \frac{d\phi}{dx}(x+\theta h)\,dx\,d\theta\\
		&=h \int_{\mathbb{R}}  \frac{d\phi}{dx}(x+\theta h)\,dx=h\,\left[\phi\right].
	\end{aligned}
\end{equation*}
Considering $h$ equal to $-c\,\textrm{dt}$ and assuming $[\phi]\neq 0$, the equality becomes
\begin{equation*}
	c=\frac{1}{\left[\phi\right]\,\textrm{dt}}\int_{\mathbb{R}} \left\{\phi(x)-\phi(x-c\, \textrm{dt})\right\}\,dx.
\end{equation*}
Assuming that $u^n_j$ is an approximation of $\phi(x_j-ct^n)$, we infer the estimate
\begin{equation}\label{numerAve}
	c\approx c^{n,1}_{\textrm{LY}}:=\frac{\mathbf{1} \cdot (u^{n} - u^{n+1})}{\left[\phi\right]}\,\cdot\,\frac{\textrm{dx}}{\textrm{dt}}
		=\frac{1}{\left[\phi\right]}\sum_{j} (u^{n}_{j} - u^{n+1}_{j})\,\cdot\,\frac{\textrm{dx}}{\textrm{dt}},
\end{equation}
where $\mathbf{1}=(1,\dots,1)$.
Hence, the value $c^n$ can be considered as a space averaged propagation speed, which
is expected to stabilize when the approximation $u^n$ converges to the given  asymptotic
profile $\phi$ with constant velocity $c$.

\subsection{Numerical experiments}
\label{subsect:numerics}

Next, we intend here to compare the results produced by the two algorithms.
In this respect, we have to specify the initial datum which will be chosen in the class of Riemann type,
i.e. corresponding to the discontinuous function
\begin{equation*}
	u_0(x)=\left\{\begin{aligned}
		1	&\quad	&x<0,\\ 0	&\quad	&x>0,\\ 
		\end{aligned}\right.
\end{equation*}
with $v_0$ determined by the corresponding values obtained by setting $\partial_t u(x,0)=0$ in the corresponding algorithm.
Such choice is very natural, since we are looking for a solution converging to the traveling front connecting the two stable state.

We focus on the case of the cubic bistable nonlinearity \eqref{polyreact} with $\alpha\in(0,1)$, 
with the goal of matching the values for the velocity $c_\ast$ as given by comparing the values provided by the exact formula
\eqref{explicitspeed2} in the case $\sigma=0$ and $\tau=1$ and the value provided by the shooting argument,
as described in Subsection \ref{subsect:comp_prop_speed}.
For sakeness of simplicity, we limit ourselves to the case $a=\kappa=1$.

We numerically solve the corresponding PDE in the space interval $[0,L]$ --with zero-flux boundary conditions-- in
the time span $[0,T]$, where we consider the case $L=50$, $T=20$ with spatial mesh $\textrm{dx}=10^{-1}$ and
time discretization $\textrm{dt}=10^{-3}$.

Finally, to quantify the error of the estimates we use the standard quantity
\begin{equation*}
	E_{\ast}:=\left|\frac{c^{n,p}_{\ast}-c_{\textrm{ex}}}{c_{\textrm{ex}}}\right|,
\end{equation*}
where $\ast\in\{\textrm{s\&s}, \textrm{LY}\}$ and $p=1$ if $\ast=\textrm{LY}$.

\subparagraph{\it Allen--Cahn equation with damping}
 Here, we compare the exact formula for the propagation speed \eqref{explicitspeed2} with the approximated estimates
obtained by applying in series one of the two scheme (first-order and Li\'enard) and, after that, the scout\&spot algorithm \eqref{levelcurvealgo}
and the LeVeque--Yee formula \eqref{numerAve}.
The results are summarized in Table \ref{tab:A}, relatively to three different choices of the intermediate (unstable) zero $\alpha$.

\begin{table}
\caption{Allen--Cahn equation with damping and polynomial reaction function, see \eqref{polyreact}:
values of $\alpha$ and $c_{\textrm{ex}}=c_{\textrm{ex}}(\alpha)$, together with the different numerical scheme,
corresponding speed estimates and relative errors.}
\label{tab:A} 
\begin{center}\begin{tabular}{p{1.00cm}p{2.00cm}p{1.50cm}p{1.00cm}p{1.75cm}p{1.50cm}p{1.50cm}p{1.50cm}}
$\alpha$	& $c_{\textrm{ex}}=c_{\textrm{ex}}(\alpha)$	& {\it scheme}	& s\&s	& $E_{\textrm{s\&s}}$	& LY			& $E_{{}_{\textrm{LY}}}$	\\
\noalign{\smallskip}\svhline\noalign{\smallskip}
0.125	&  0.4685213		& first-order	& 0.47	& $3.16\times 10^{-3}$	& 0.4682076	& $6.69\times 10^{-4}$	\\
		&				& Li\'enard	& 0.46	& $1.82\times 10^{-2}$	& 0.4662342	& $4.88\times 10^{-3}$	\\
\noalign{\smallskip}\hline\noalign{\smallskip}
0.250	& 0.3333333		&first-order	& 0.34	& $2.00\times 10^{-2}$	& 0.3331151	& $6.55\times 10^{-4}$	\\
		&				& Li\'enard	& 0.33	& $1.00\times 10^{-2}$	& 0.3310495	& $6.85\times 10^{-3}$	\\
\noalign{\smallskip}\hline\noalign{\smallskip}
0.375	& 0.1740777		&first-order	& 0.17	& $2.34\times 10^{-2}$	& 0.1739747	& $5.92\times 10^{-4}$	\\
		&				& Li\'enard	& 0.17	& $2.34\times 10^{-2}$	& 0.1715496	& $1.45\times 10^{-3}$	\\
\noalign{\smallskip}\hline\noalign{\smallskip}
\end{tabular}\end{center}
\end{table}

It is transparent the higher precision of the LeVeque--Yee formula \eqref{numerAve} which add to the number of free parameters to be chosen
in the scout\&spot algorithm (such as the level $\theta$, the value of $p$...), making the use of the latter strategy less effective.

Next, we pass to analyze the Allen--Cahn equation with a piecewise linear reaction function
with a jump point located at $u=\alpha$.
In this case, the crucial problem is, of course, the presence of a discontinuity in the source term.
Thus, we compare the capability of the scout\&spot algorithm and the LeVeque--Yee formula.
The results, obtained by using the same numerical data previously described, are reported in Table \ref{tab:Apwl}.
As can be appreciated from the values, the error is always of the order of $1\%$, which is largely acceptable.
\begin{table}[htb]
\caption{Allen--Cahn equation with damping and piecewise affine reaction function, see \eqref{pwl}:
Values of $\alpha$ and $c_{\textrm{ex}}=c_{\textrm{ex}}(\alpha)$, together with the different numerical scheme,
corresponding speed estimates and relative errors.}
\label{tab:Apwl} 
\begin{center}\begin{tabular}{p{1.00cm}p{2.00cm}p{1.50cm}p{1.00cm}p{1.75cm}p{1.50cm}p{1.50cm}p{1.50cm}}
$\alpha$	& $c_{\textrm{ex}}=c_{\textrm{ex}}(\alpha)$	& {\it scheme}	& s\&s	& $E_{\textrm{s\&s}}$	& LY		& $E_{{}_{\textrm{LY}}}$	\\
\noalign{\smallskip}\svhline\noalign{\smallskip}
0.125	& 0.9149914		& first-order	& 0.90	& $1.64\times 10^{-2}$	& 0.9021793	& $1.40\times 10^{-2}$	\\
		&				& Li\'enard	& 0.90	& $1.64\times 10^{-2}$	& 0.9006799	& $1.56\times 10^{-2}$	\\
\noalign{\smallskip}\hline\noalign{\smallskip}
0.250	& 0.7559289		& first-order	& 0.74	& $2.11\times 10^{-2}$	& 0.7496325	& $8.33\times 10^{-3}$	\\
		&				& Li\'enard	& 0.74	& $2.11\times 10^{-2}$	& 0.7484820	& $9.85\times 10^{-3}$	\\
\noalign{\smallskip}\hline\noalign{\smallskip}
0.375	& 0.4588315		& first-order	& 0.45	& $1.92\times 10^{-2}$	& 0.4557922	& $6.62\times 10^{-3}$	\\
		&				& Li\'enard	& 0.46	& $2.55\times 10^{-3}$	& 0.4554450	& $7.38\times 10^{-3}$	\\
\noalign{\smallskip}\hline\noalign{\smallskip}
\end{tabular}\end{center}
\end{table}

As shown by the numerical results, also the case of a discontinuous reaction function can be handled by both algorithms,
with  slightly better error estimates for the LeVeque--Yee formula (which is also very easy to implement).

\subparagraph{Allen--Cahn equation with relaxation}
Finally, we consider the case of the hyperbolic Allen--Cahn equation with relaxation, that is \eqref{scalar_damp_rd_1d}
with $\sigma=\tau>0$ (fixed equal to 1, for simplicity) for the third order polynomial reaction function, given by \eqref{polyreact}.
In such a case, in addition to the first-order and Li\'enard schemes, we may also apply the kinetic scheme, also presented in Subsection \ref{subsec:algorithms}.
A selection of the results are collected in Table \ref{tab:B} and confirm the same conclusion as above: 
with the same space-time grid, the LeVeque--Yee formula is to be preferred, since it guarantees greater precision in speed approximation.

\begin{table}
\caption{Values of $\alpha$ and $c_{\textrm{ex}}=c_{\textrm{ex}}(\alpha)$, together with the different numerical schemes,
speed estimates and relative errors.}
\label{tab:B} 
\begin{center}\begin{tabular}{p{1.00cm}p{2.00cm}p{1.50cm}p{1.00cm}p{1.75cm}p{1.50cm}p{1.50cm}p{1.50cm}}
$\alpha$	& $c_{\textrm{ex}}=c_{\textrm{ex}}(\alpha)$ & {\it scheme}	& s\&s	& $E_{\textrm{s\&s}}$	& LY			& $E_{{}_{\textrm{LY}}}$	\\
\noalign{\smallskip}\svhline\noalign{\smallskip}
0.125	&  0.5342843	& first-order	& 0.53	& $8.02\times 10^{-3}$	& 0.5335445	& $1.38\times 10^{-3}$	\\
		&			& Li\'enard	& 0.53	& $8.02\times 10^{-3}$	& 0.5318317	& $4.59\times 10^{-3}$	\\
		&			& kinetic		& 0.54	& $1.07\times 10^{-2}$	& 0.5347508	& $8.73\times 10^{-4}$	\\
\noalign{\smallskip}\hline\noalign{\smallskip}
0.250	& 0.3754283	& first-order	& 0.38	& $1.22\times 10^{-2}$	& 0.3750573	& $9.88\times 10^{-4}$	\\
		&			& Li\'enard	& 0.37	& $1.45\times 10^{-2}$	& 0.3728276	& $6.93\times 10^{-3}$	\\
		&			& kinetic		& 0.38	& $1.22\times 10^{-2}$	& 0.3758528	& $1.13\times 10^{-3}$	\\
\noalign{\smallskip}\hline\noalign{\smallskip}
0.375	& 0.1941490	& first-order	& 0.19	& $2.14\times 10^{-2}$	& 0.1940086	& $7.23\times 10^{-4}$	\\
		&			& Li\'enard	& 0.19	& $2.14\times 10^{-2}$	& 0.1913620	& $1.44\times 10^{-3}$	\\
		&			& kinetic		& 0.19	& $2.14\times 10^{-2}$	& 0.1943773	& $1.18\times 10^{-3}$	\\
\noalign{\smallskip}\hline\noalign{\smallskip}
\end{tabular}\end{center}
\end{table}

Other numerical experiments have been performed with different choices of $p$ and better precision for the estimate of the scout\&spot algorithm,
providing a corrisponding higher order of precision of the LeVeque--Yee formula, which appear again as a more precise tool.
Comparing the three types of scheme --first-order reduction, Li\'enard, kinetic-- the first two have some very poor resolution of the equation for short time,
in particular when considered in relation with the third one.
Spurious oscillations are generated by both the schemes due to the presence of a discontinuity in the initial datum.
Differently, the kinetic algorithm is capable of reproducing the correct behavior also in the short time (see \cite{LattMascPlazSime16,LattMascPlazSime19} 
for more numerical simulations).
Nevertheless, we stress that the latter is much slower with respect to the other two.
Thus, computing the propagation speed --which is a parameter relevant for the large-time behavior-- the short time behavior is of secondary importance
with respect to the capability of the scheme of being capable to reproduce the main features of the model in the long run, once the evolution has already 
solved the initial problem of the presence of a jump.
This is particularly crucial because of the presence of the reaction term which, in large part of the space, pushes the solutions to stay close to stable solution of
the underlying ODE.

The case of the piecewise affine reaction function, described in the last paragraph of Subsection \ref{subsect:specialcases}, is harder to be simulated,
since the numerical schemes of Subsection \ref{subsec:algorithms} are not well-behaved in the presence of discontinuous reaction function due to the presence
of the term $\tau f(u)$ differentiated with respect to time.
Numerical deficiencies arise already when performing simulations of the PDE, inherited by the jump of the reaction function $f$, probably due to the
fact that the dissipativity condition $1-\tau f'>0$ is never satisfied at $\alpha$ whenever $\tau>0$,
At the moment, we are not aware of any numerical schemes which is capable of performing reliable simulations in presence of discontinuities.

\begin{acknowledgement}
The authors are thankful to the anonymous referee for a number of significant remarks which drastically improved the content of the manuscript.
Simulations have been performed by {\sc Scilab 6.0.2}, \url{https://www.scilab.org/}.
\end{acknowledgement}

\end{document}